\documentclass[centertags,leqno]{article}

\usepackage{amsmath}
\usepackage{amssymb}
\usepackage{latexsym}
\usepackage{amsxtra}
\usepackage{amscd}
\usepackage{theorem}
\usepackage{epic}
\usepackage{eepic}
\usepackage[all]{xy}

\theoremstyle{plain}

{\theorembodyfont{\slshape}

        \newtheorem{thm}{Theorem}[section]
        \newtheorem{cor}[thm]{Corollary}
        \newtheorem{lem}[thm]{Lemma}
        \newtheorem{prop}[thm]{Proposition}
}

{\theorembodyfont{\rmfamily}

        \newtheorem{defn}[thm]{Definition}
        \newtheorem{prob}[thm]{Problem}
        \newtheorem{rem}[thm]{Remark}
        
        \newtheorem{ass}[thm]{Assumption}
        \newtheorem{notation}[thm]{Notation}

        \newtheorem{constr}[thm]{Construction}
        
        \newtheorem{result}[thm]{Result}
}

\renewcommand{\em}{\sl}

\makeatletter
\renewcommand{\subsubsection}{\@startsection{subsubsection}{3}%
        {\z@}{-3.25ex plus -1ex minus-.2ex}{-1em}{\bf}}
\renewcommand{\section}{\@startsection{section}{1}%
        {\z@}{-\baselineskip}{\baselineskip}{\bf\large}}
\renewcommand{\subsection}{\@startsection{subsection}{2}%
        {\z@}{-\baselineskip}{\baselineskip}{\bf}}
\makeatother

\newcommand{\proof}{{\bf Proof:\ }}
\newcommand{\Endproof}{\hspace*{\fill} $\Box$ \vspace{1ex} \noindent }

\newcommand{\ZZ}{\mathbb{Z}}

\newcommand{\PP}{\mathbb{P}}

\newcommand{\FF}{\mathbb{F}}

\newcommand{\HH}{\mathbb{H}}

\newcommand{\RR}{\mathbb{R}}
\newcommand{\C}{\mathfrak{C}}

\newcommand{\Oo}{\mathcal{O}}
\newcommand{\A}{\mathcal{A}}
\newcommand{\B}{\mathcal{B}}

\newcommand{\D}{\mathcal{D}}
\newcommand{\U}{\mathcal{U}}

\newcommand{\E}{\mathcal{E}}

\newcommand{\Od}{\hat{\Oo}}

\newcommand{\G}{\mathcal{G}}
\newcommand{\F}{\mathcal{F}}

\newcommand{\Ll}{\mathcal{L}}

\newcommand{\K}{\mathcal{K}}

\newcommand{\I}{\mathcal{I}}
\newcommand{\J}{\mathcal{J}}

\newcommand{\T}{\mathcal T}
\newcommand{\M}{\mathcal{M}}

\newcommand{\X}{\mathcal{X}}
\newcommand{\Y}{\mathcal{Y}}

\newcommand{\Fu}{\underline{F}}

\newcommand{\At}{\tilde{A}}

\newcommand{\Spec}{\mathop{\rm Spec}} 
\newcommand{\Hom}{\mathop{\rm Hom}\nolimits} 
\newcommand{\HOm}{\mathop{\mathcal{H}om}\nolimits} 
\newcommand{\RHom}{\mathop{\RR \rm Hom}\nolimits} 
\newcommand{\RHOm}{\mathop{\RR \mathcal{H}om}\nolimits} 
\newcommand{\HOM}{\mathop{\mathfrak{Hom}}\nolimits}

\newcommand{\HExt}{\mathop{\rm \mathbb{E}xt}\nolimits} 
\newcommand{\EXT}{\mathop{\mathfrak{Ext}}\nolimits} 
\newcommand{\EEXT}{\underline{\EXT}}

\newcommand{\Mod}{\mathfrak{M}} 
\newcommand{\Der}{\mathfrak{D}} 
\newcommand{\Ko}{\mathfrak{K}} 
\newcommand{\Ab}{\mathfrak{Ab}} 
\newcommand{\Tot}{\mathrm{Tot}} 
\newcommand{\EXt}{\mathop{\mathcal{E}xt}\nolimits}

\newcommand{\Id}{\mathop{\rm Id}\nolimits} 
\newcommand{\Pic}{\mathfrak{P}}
\newcommand{\PPic}{\underline{\Pic}}
\newcommand{\PIC}{\mathop{\mathfrak{Pic}}\nolimits} 
\newcommand{\PPIC}{\underline{\PIC}}

\newcommand{\ord}{\mathop{\rm ord}\nolimits}

\renewcommand{\Im}{\mathop{\rm Im}} 
\newcommand{\Ker}{\mathop{\rm Ker}} 
\newcommand{\Def}{\mathop{\rm Def}\,}

\newcommand{\univ}{_{\rm\scriptscriptstyle univ}}

\newcommand{\prim}{_{\rm\scriptscriptstyle prim}}
\newcommand{\new}{_{\rm\scriptscriptstyle new}}

\newcommand{\tame}{_{\rm\scriptscriptstyle tame}}
\newcommand{\wild}{_{\rm\scriptscriptstyle wild}}
\newcommand{\branch}{_{\rm\scriptscriptstyle ram}}
\renewcommand{\int}{_{\rm\scriptscriptstyle int}}

\newcommand{\diff}{{\rm d}}

\newcommand{\Rt}{\tilde{R}}

\newcommand{\Xd}{\hat{X}}
\newcommand{\Yd}{\hat{Y}}

\newcommand{\Ad}{\hat{A}}
\newcommand{\Bd}{\hat{B}}
\newcommand{\Ld}{\hat{L}}
\newcommand{\Ih}{\hat{I}}
\newcommand{\Mh}{\hat{M}}

\newcommand{\ub}{\bar{u}}
\newcommand{\cb}{\bar{c}}

\newcommand{\phib}{\bar{\phi}}

\newcommand{\m}{\mathfrak{m}}
\renewcommand{\a}{\mathfrak{a}}

\newcommand{\inj}{\hookrightarrow}
\newcommand{\To}{\;\longrightarrow\;}
\newcommand{\iso}{\stackrel{\sim}{\to}}
\newcommand{\liso}{\;\stackrel{\sim}{\longrightarrow}\;}

\newcommand{\pfeil}[1]{\stackrel{#1}{\to}}
\newcommand{\lpfeil}[1]{\stackrel{#1}{\To}}

\newcommand{\bmu}{\boldsymbol{\mu}}

\newcommand{\bi}{\underline{i}}

\newcommand{\gen}[1]{\mathopen\langle#1\mathclose\rangle}
\newcommand{\vf}[1]{{\textstyle \frac{\diff}{\diff #1}}}

\newcommand{\rupfeil}[1]{\Big\downarrow\vcenter{%
                         \rlap{$\scriptstyle #1$}}}

\newcommand{\Prod}{\displaystyle\prod}

\newcommand{\dirlim}[1]{\underset{\underset{#1}{\To}}{\lim}}

\makeatletter
\renewcommand{\subsection}{\@startsection{subsection}{2}%
        {\z@}{-3.25ex plus -1ex minus-.2ex}{-1em}{\bf}}
\makeatother

\title{Formal deformation of curves with group scheme action}
\author{Stefan Wewers\\[1ex]
        Universit\"at Bonn}
\date{}

\begin{document}

\maketitle

\begin{abstract}
  We study equivariant deformations of singular curves with an action
  of a finite flat group scheme, using a simplified version of
  Illusie's equivariant cotangent complex. We apply these methods in a
  special case which is relevant for the study of the stable reduction
  of three point covers.
\end{abstract}


\section*{Introduction}

Let $Y$ be a (not necessarily smooth) projective curve over an
algebraically closed field $k$ of characteristic $p>0$. Let
$W$ be a complete local ring with residue field $k$.
Furthermore, let $G$ be a finite flat group scheme over $W$
which acts faithfully on $Y$. We denote by $\Def(Y,G)$ the
functor which associates to a local Artinian $W$-algebra $R$
with residue field $k$ the set of isomorphism classes of
$G$-equivariant deformations of $Y$ to $R$. The present
paper is concerned with a study of the functor $\Def(Y,G)$,
using cohomological methods. The special case where $Y$ is smooth
and $G$ is a constant group scheme has been studied
previously by Bertin and M\'ezard \cite{BertinMezard00}.

One of the motivations for studying the functor $\Def(Y,G)$
is the {\em lifting problem}. Suppose that $Y$ is smooth and
that $G$ is a finite abstract group which acts faithfully on
the curve $Y$. Let $W$ be the ring of Witt vectors over $k$,
and consider $G$ as a constant group scheme over $W$. In
this situation, the lifting problem ask the following
question. Does there exist a finite extension $R/W$ of
complete discrete valuation rings and a $G$-equivariant lift
of $Y$ over $R$? For instance, if the deformation functor
$\Def(Y,G)$ is unobstructed then the answer to this question
is positive.

A conjecture of Oort predicts that the lifting problem has a positive
solution if the group $G$ is cyclic.  However, even in the simplest
nontrivial case $G=\ZZ/p$ (where Oort's conjecture is proved, see
\cite{SOS} and \cite{GreenMatignon98}) the functor $\Def(Y,G)$ is
obstructed. In \cite{BertinMezard00} these obstructions are identified
as elements in a certain Galois cohomology group. However, they remain
a bit mysterious. One of the motivations for generalizing the approach of
Bertin--M\'ezard is the author's hope that this will lead to new
insight into the nature of these obstruction, and the lifting problem
in general. 

Another (related) motivation comes from the study of the
stable reduction of Galois covers of curves. Let $R$ be a
complete discrete valuation ring, with algebraically closed
residue field $k$ of characteristic $p$ and fraction field
$K$ of characteristic $0$. Let $Y_K\to X_K$ be a Galois
cover of smooth projective curves over $K$, with Galois
group $G$.  After a finite extension of $K$, there exists a
certain natural $R$-model $Y_R\to X_R$ of $Y_K\to X_K$,
called the {\em stable model}, see \cite{Raynaud98} or
\cite{bad} for a precise definition. The problem we are
interested in is to understand this model and in particular
its relation with the ramification of the prime $p$ in the
field $K$.  It has become clear from recent work of several
authors (e.g.\ \cite{YannickArbres}, \cite{SaidiTorsors},
\cite{bad}) that this problem naturally leads to the study
of singular curves with an action of an infinitesimal group
scheme, and of the deformation theory of such objects.

\vspace{1ex} This paper is divided into two main parts.  The first
part (\S 1--3) is an exposition of certain cohomological methods for
studying equivariant deformations of (not necessarily smooth) curves
with group scheme action. Although the guiding principles are the same
as in \cite{BertinMezard00}, we have to use much heavier technical
machinery.  For instance, it does not suffice to look at the
equivariant cohomology of the tangent bundle of $Y$ over $k$, as in
\cite{BertinMezard00}. Instead, one has to consider certain hyperext
groups with values in the {\em equivariant cotangent complex} of $Y$
over $k$. The latter is an object in the derived category of
$G$-$\Oo_Y$-modules, and was first introduced by Grothendieck in
\cite{GrothendieckCC}.

In principal, everything one might want to known about the
equivariant cotangent complex and its role in deformation theory can
be found in Illusie's book \cite{IllusieCC}. However, the generality
in which \cite{IllusieCC} is written makes it difficult to
read and to work with in a concrete situation (at least for the author
of this paper). In the literature there are a number of excellent and
readable accounts of certain special cases (see e.g.\ \cite{Vistoli99}
or \cite{BertinMezard00}, \S 2-3) but none seems to be sufficiently
general to deal with the case we need.  To improve this situation, the
present paper contains a self-contained exposition of a special case
of Illusie's theory, which should nevertheless be sufficiently general
for the applications we have in mind.

\vspace{1ex} In the second part of this paper (\S 4--5) we apply the
general theory to a special case which is relevant for the study of
three point covers with bad reduction. In particular, we prove a
certain result which is a key ingredient for the main theorem of
\cite{bad}.

We start in \S 4 with a {\em multiplicative deformation datum}. To
give an idea what this is, fix a smooth projective curve $X$ over $k$.
Then a multiplicative defirmation datum over $X$ is a pair $(Z,V)$,
where $Z\to X$ is a Galois cover of smooth projective curves over a
field $k$ of characteristic $p>0$, with Galois group $H$ of order
prime-to-$p$, and $V\subset\Omega_{k(Z)/k}$ is an $H$-stable
$\FF_p$-vector space of logarithmic differential forms. To $(Z,V)$ we
associate a finite flat group scheme $G$ over $W(k)$ and a (singular)
curve $Y$ over $k$ with an action of $G$ such that $X=Y/G$.  Briefly,
the group scheme $G$ is of the form $\bmu_p^s\rtimes H$ and $Y\to Z$
is the $\bmu_p^s$-cover locally given by $s$ Kummer equations
$y_i^p=u_i$, where $\phi_i=\diff u_i/u_i$, $i=1,\ldots,s$, form a
basis of $V$.  

We study the deformation functor $\Def(Y,G)$ which
classifies equivariant formal deformations of $Y$ and exhibit a number
of its properties which are, in general, very different from the
properties enjoyed by the deformation functor studied in
\cite{BertinMezard00}. For instance, there is in general no such thing
as a local-global principle, because the `local contribution' to the
tangent space of the functor $\Def(Y,G)$ is not concentrated in a
finite number of closed points. However, from another point of view
things are really much easier than in \cite{BertinMezard00}, due to
the fact that the `$p$-Sylow' of $G$ is a multiplicative group scheme.
Since multiplicative group schemes have trivial cohomology, the
general machinery developed in the first sections shows that the
deformation functor $\Def(Y,G)$ is unobstructed.  Another nice
property of $\Def(Y,G)$ is the existence of a natural morphism of
deformation functors
\begin{equation} \label{introeq1}
     \Def(Y,G) \;\To\; \Def(X;\tau_j)
\end{equation}
which sends an equivariant deformation of $Y$ to its quotient by $G$.
(Here we regard $X$ as a marked curve, the marked points being the
`branch points' $\tau_1,\ldots,\tau_n$ of the $G$-cover $Y\to X$.)
In this respect, the $G$-cover $Y\to X$ behaves like a tamely ramified
Galois cover. However, unlike in the case of tamely ramified Galois
covers, the functor \eqref{introeq1} is in general not an isomorphism.

\vspace{1ex} In \S 5 we assume in addition that the curve $X$ is the
projective line and that the vector space $V$ is an irreducible
$\FF_p[H]$-module which decomposes, after tensoring with
$\bar{\FF}_p$, into the direct sum of one dimensional modules. Among
all the multiplicative deformation data $(Z,V)$ of this type, there
are some which we call {\em special}.  The definition of specialty is
given in terms of certain numerical invariants attached to $(Z,V)$.
But philosophically, special deformation data are attached to three
point Galois covers of the projective line with bad reduction to
characteristic $p$. We refer to \cite{special} and \cite{bad} for
details on the case $\dim_{\FF_p}V=1$ and for a more satisfactory
explanation of the connection to three point covers. Let us only
mention that the deformation theory of the $G$-cover $Y\to X$ attached
to a special deformation datum has a number of very nice and
surprising properties:
\begin{itemize}
\item
  {\em The lifting property:} the morphism of deformation functors
  \eqref{introeq1} is an isomorphism. In this respect, the $G$-cover
  $Y\to X$ behaves just like a tamely ramified Galois cover. 
\item {\em The local-global principle:} local deformations in formal
  neighborhoods of the ramification points (which satisfy a certain condition)
  can be interpolated by a {\em unique} global deformation of $Y$. 
\item
  {\em Rigidity:} If an equivariant deformation of $Y$ in equal
  characteristic (i.e.\ over a local $k$-algebra) is again special 
  then it is the trivial deformation.
  Therefore, there exist at most a finite number of special
  deformation data of a given type (up to isomorphism), and every
  special deformation datum can be defined over a finite field.
\end{itemize}
These properties are very particular to special deformation data. They
reflect, in a rather subtle way, the connection to three point covers
with bad reduction and in particular to the fact that three point
covers are `rigid' objects.

\vspace{1ex}
At the end of the paper, the reader will find three appendices
containing background material which the author found difficult to
extract from the literature. This includes Picard stacks, the
cohomology of affine group schemes, and two spectral sequences which
are useful to compute equivariant hyperext groups.


\section{The equivariant cotangent complex} \label{constr}

In \cite{IllusieCC} Illusie defines, for any morphism of schemes $Y\to
S$, the cotangent complex $\Ll_{Y/S}$. This is a complex of flat
$\Oo_Y$-modules, well defined up to canonical quasi-isomorphism, such
that $H^0(\Ll_{Y/S})=\Omega_{Y/S}$. If $Y\to S$ is smooth then
$\Ll_{Y/S}=\Omega_{Y/S}$. Moreover, if $G\to S$ is a group scheme
acting on $Y$, Illusie defines the {\em equivariant cotangent complex}
as an object of the derived classifying topos $\Der^+(BG_{/X})$ whose
underlying complex of $\Oo_Y$-modules is $\Ll_{Y/S}$. 

In this section we give a more down-to-earth definition of the
equivariant cotangent complex which, however, works well only if $Y\to
S$ and the $G$-action on $Y$ have certain good properties. We follow
the original approach of Grothendieck \cite{GrothendieckCC}.  This
gives the `correct' cotangent complex only if $Y\to S$ is a local
complete intersection morphism. We assume that $Y$ admits locally an
equivariant embedding into a formally smooth $S$-scheme with
$G$-action. Under this assumption, it is much easier to endow the
cotangent complex with a natural $G$-action.

\subsection{} \label{constr1}

Let $S=\Spec R$ be an affine scheme, $G\to S$ a flat affine group
scheme and $Y\to S$ an $S$-scheme with an action of $G$.  By a {\em
  $G$-$\Oo_Y$-module} we mean a sheaf of $\Oo_Y$-modules $\F$, together
with a lift of the $G$-action from $Y$ to $\F$.  A homomorphism between
two $G$-$\Oo_Y$-modules $\F$ and $\G$ is a sheaf homomorphism which is
both $\Oo_Y$-linear and $G$-equivariant. The group of such
homomorphisms is denoted by $\Hom_G(\F,\G)$. We denote by $\Mod(Y,G)$
the corresponding category of $G$-$\Oo_Y$-modules.  See Appendix
\ref{GOY2} for more details on the category $\Mod(Y,G)$. For
$*\in\{+,-,b\}$, we denote by $\Ko^*(Y,G)$ the category of cochain
complexes in $\Mod(Y,G)$, up to homotopy, which are bounded from below
($*=+$), bounded from above ($*=-$) or bounded in both directions
($*=b$). We write $\Der^*(Y,G)$ for the derived category of
$\Ko^*(Y,G)$.

In this section we define the equivariant cotangent complex
$\Ll_{Y/S}$ of the morphism $Y\to S$ as an object of $\Der^+(Y,G)$,
assuming:

\begin{ass} \label{mainass}
  Every point of $Y$ is contained in an affine and $G$-stable open
  neighborhood $U\subset Y$ such that the following holds. There
  exists a formally smooth affine $S$-scheme $P\to S$ with $G$-action
  and a $G$-equivariant closed immersion $\varphi:U\inj P$.
\end{ass}

\begin{rem}
  It is not clear to the author how restrictive Assumption
  \ref{mainass} is. We expect that it can be verified in any concrete
  situation where one actually wants to apply our theory. For
  instance, in \S \ref{defdat} we use the case where $G$ is an
  extension of a constant by a multiplicative group scheme and acts
  freely on a dense open subset of $Y$. In this situation, Assumption
  \ref{mainass} is easy to verify.
\end{rem}


\subsection{}   \label{localconstr}

A triple $(U,P,\varphi)$ as in Assumption \ref{mainass} is called a
{\em local chart} for $Y\to S$. Often we will simply write $\varphi$
instead of $(U,P,\varphi)$. Given such a local chart, we denote by
$\I\subset\Oo_P$ the sheaf of ideals defining the image of $\varphi$.
We define the {\em cotangent complex} of the chart $\varphi$ as the
following complex of $G$-$\Oo_Y$-modules:
\begin{equation} \label{cotangeq1}
    \Ll_{\varphi} \;:=\; (\I/\I^2 \To \Omega_{P/S}\otimes\Oo_Y).
\end{equation}
The two nontrivial terms of $\Ll_\varphi$ lie in degree $-1$ and $0$.
Note that there is a natural augmentation $\Ll_{\varphi}\to
\Omega_{Y/S}$ which identifies $\Omega_{Y/S}$ with
$H^0(\Ll_{\varphi})$. 

\begin{rem}
  \begin{enumerate}
  \item 
    If $Y/S$ is of finite type, then we may take $P/S$ to be
    smooth. In this case, $\Ll_\varphi^0=\Omega_{P/S}\otimes\Oo_Y$ is a
    locally free $\Oo_Y$-module of finite rank.
  \item If, moreover, $Y\to S$ is a local complete intersection (in
    the sense of \cite{SGA6}, VIII.1.1) then the embedding $\varphi$
    is regular. Recall that this means the following. For every point
    $y\in U$ the stalk $\I_y$ is an ideal generated by a regular
    sequence of the local ring $\Oo_{P,y}$. It follows that
    $\Ll_\varphi^{-1}=\I/\I^2$ is a locally free $\Oo_Y$-module of finite
    rank, too. 
  \end{enumerate}
\end{rem}

Let $(U,P,\varphi)$ and $(V,Q,\psi)$ be two local charts, and assume
that $V\subset U$. A {\em morphism} from $\psi$ to $\varphi$ is a
$G$-equivariant morphism of $S$-schemes $u:Q\to P$ such that the
diagram
\[\begin{CD}
   V  @>{\psi}>>     Q        \\
   @VVV     @VV{u}V  \\
   U  @>{\varphi}>>  P        \\
\end{CD}\]
commutes. We use the notation $u:\psi\to\varphi$.  Note that $u$
induces a morphism of complexes of $G$-$\Oo_Y$-modules
\[
       u^*:\;\Ll_{\varphi}|_V \;\To\; \Ll_\psi.
\]

\begin{lem} \label{constr2lem}
\begin{enumerate}
\item
  The homotopy class of $u^*$ is independent of the morphism $u$.
\item
  The morphism $u^*$ is a quasi-isomorphism.
\end{enumerate}
\end{lem}

\proof It is no restriction to assume that $U=V$. We may also assume
that $U$ is affine. Let $P'$ denote the second infinitesimal
neighborhood of $U$ in $P$, i.e.\ the closed subscheme of $P$ defined
by the sheaf of ideals $\I^2$. Similarly, $Q'$ denotes the second
infinitesimal neighborhood of $Y$ in $Q$. It is defined by $\J^2$,
where $\J\subset\Oo_Q$ is the sheaf of ideals defining the image of
$\psi$. Let $v:Q\to P$ be another morphism of local charts, and set
$u':=u|_{Q'}$ and $v':=v|_{Q'}$. It is clear that $u^*$ (resp.\ $v^*$)
only depends on the restriction $u'$ (resp.\ on $v'$). An easy
computation shows that the difference of the two pullback maps
\[
       (u')^*-(v')^*:\;\Oo_P \;\To\; \J/\J^2
\]
is an $R$-linear derivation. Hence it gives rise to
an $\Oo_Y$-linear map $s:\Omega_{P}\otimes\Oo_Y\to\J/\J^2$, and one
checks that $s$ is the desired homotopy between $u^*$ and $v^*$:
\[
  \xymatrix{
     \I/\I^2 \ar[r]^-d \ar[d]_{u^*-v^*}
           & \Omega_{P/S}\otimes\Oo_Y  \ar[ld]_s
                \ar[d]^{u^*-v^*}  \\ 
     \J/\J^2 \ar[r]^-d   & \Omega_{Q/S}\otimes\Oo_Y }
\]
This proves (i).

By assumption $Q\to S$ is formally smooth and $U$ is affine. Hence
there exists a morphism $w':P'\to Q$ lifting $\psi:Y\inj Q$. As in the
proof of (i), one shows that there are homotopies
\[
         (w')^*\circ u^* \;\sim\; \Id_{\Ll_\varphi}, \qquad
         u^*\circ(w')^* \;\sim\; \Id_{\Ll_{\psi}}.
\]
This proves (ii).
\Endproof

\subsection{} \label{globalconstr}

We are now ready to define the equivariant cotangent complex. 
By Assumption \ref{mainass} there exists a covering $(U_i)_{i\in
  I}$ of $Y$ by affine and $G$-stable opens $U_i\subset Y$, each
admitting a local chart $\varphi_i:U_i\inj P_i$.  We choose, once and
for all, a well-ordering on the set of indices of the covering
$(U_i)$.  The datum $(U_i,\varphi_i)$ is called an {\em atlas}.

For any $(n+1)$-tuple $\bi=(i_0,\ldots,i_n)$ we set
\[
     U_{\bi}:=U_{i_0}\cap\ldots\cap U_{i_n}, \qquad
     P_{\bi}:=P_{i_0}\times_S\ldots\times_SP_{i_n}, \qquad
     \varphi_{\bi}:=\varphi_{i_0}\times\cdots\times\varphi_{i_n}:\,
        U_{\bi}\inj P_{\bi}.
\]
Clearly, $\varphi_{\bi}$ is a local chart and gives rise to a complex
of $G$-$\Oo_{U_{\bi}}$-modules $\Ll_{\varphi_{\bi}}$. We denote by
$\Ll_{\bi}$ the push-forward of $\Ll_{\varphi_{\bi}}$ to $Y$. Thus,
$\Ll_{\bi}$ is a flat and quasi-coherent $G$-$\Oo_Y$-module such that
$\Ll_{\bi}|_{U_i}=\Ll_{\varphi_{\bi}}$.

For $\bi=(i_0,\ldots,i_n)$ as above and $0\leq\nu\leq n$, let
$p_{\bi}^\nu:P_{\bi}\to P_{\bi'}$ denote the projection which leaves
out the $\nu$th component (i.e.\ 
$\bi'=(i_0,\ldots,\widehat{i_\nu},\ldots,i_n)$); it is a morphism
$\varphi_{\bi}\to\varphi_{\bi'}$ of local charts. The
resulting morphism
$(p_{\bi}^\nu)^*:\Ll_{\varphi_{\bi'}}|_{U_{\bi}}\to\Ll_{\varphi_{\bi}}$
extends in a canonical way to a morphism of $G$-$\Oo_Y$-modules
$\partial_{\bi}^\nu:\Ll_{\bi'}\to \Ll_{\bi}$.  Note that
\begin{equation} \label{localconstreq1}
  \partial_{\bi}^{\nu}\circ\partial_{\bi'}^\mu \;=\;
  \partial_{\bi}^{\mu}\circ\partial_{\bi''}^{\nu-1}
\end{equation}
holds for $\mu<\nu$, if we set
$\bi':=(\ldots,\widehat{i_\nu},\ldots)$ and
$\bi'':=(\ldots,\widehat{i_\mu},\ldots)$. 

\begin{defn} \label{Ldef}
  The {\em equivariant cotangent complex} of the morphism $Y\to S$
  (relative to the open covering $(U_i)$ and the local charts
  $\varphi_i$) is the total complex
  \[
     \Ll_{Y/S} \;:=\; \Tot(\K)
  \]
  of the following double complex of $G$-$\Oo_Y$-modules:
  \[
       \K:\quad\left\{\quad
  \begin{array}{ccc} 
    \Prod_i\;\Ll_i^{-1} & \lpfeil{d} & \Prod_i\;\Ll_i^0  \\ 
    \rupfeil{\partial}  & &  \rupfeil{\partial}        \\ 
     \Prod_{i<j}\;\Ll_{i,j}^{-1} & \lpfeil{-d} & 
       \Prod_{i<j}\;\Ll_{i,j}^0  \\
    \rupfeil{\partial}  & &  \rupfeil{\partial}  \\ 
     \Prod_{i<j<k}\;\Ll_{i,j,k}^{-1} & \lpfeil{d} &  
              \Prod_{i<j<k}\;\Ll_{i,j,k}^0 \\
    \rupfeil{\partial}  & &  \rupfeil{\partial} \\
      \cdots  & & \cdots 
  \end{array}
    \right.  
  \]
  The vertical differentials are defined as
  $\partial:=\sum_{\nu=0}^p(-1)^\nu
  \prod_{\bi}\partial_{\bi}^\nu:\K^{p,q}\to \K^{p+1,q}$. The
  horizontal differentials are induced from the differentials of the
  complexes $\Ll_{\bi}$. The columns of $\K$ start with degree $0$, so
  $\Ll_{Y/S}$ starts with degree $-1$. Note that $\Ll_{Y/S}$
  consists of flat and quasi-coherent sheaves.
\end{defn}

\begin{prop} \label{constr2prop}
  For all $i$ there exist a quasi-isomorphism
  $\beta_i:\Ll_{Y/S}|_{U_i}\to\Ll_{\varphi_i}$. Moreover, for all
  $i<j$ we have a homotopy
  \[
    s_{i,j}:\;\partial_{i,j}^1\circ\beta_i \;\sim\; 
         \partial_{i,j}^0\circ\beta_j.
  \]
  which satisfies the cocycles relation
  \begin{equation} \label{cocyclerel1}
      \partial_{i,j,k}^0\circ s_{j,k} - 
      \partial_{i,j,k}^1\circ s_{i,k} +
      \partial_{i,j,k}^2\circ s_{i,j} \;=\; 0.
  \end{equation}
\end{prop}

\proof The natural projections $\K^{0,q}|_{U_i}\to\Ll_i^q$ induce a
morphism $\beta_i:\Ll_{Y/S}|_{U_i}\to\Ll_{\varphi_i}$. 
We define the homotopy $s_{i,j}$ as follows:
\[
  s_{i,j}^0:\;\left\{
  \begin{array}{ccc}
      (\Ll_{Y/S})^0 & \;\To\; & \Ll_{i,j}^{-1} \\
      (f_k;g_{k,l})  & \;\longmapsto\; & g_{i,j}   \\
  \end{array}\right., \qquad
  s_{i,j}^1:\;\left\{
  \begin{array}{ccc}
      (\Ll_{Y/S})^1 & \;\To\; & \Ll_{i,j}^{0}        \\
      (f_{k,l};g_{k,l,m})  & \;\longmapsto\; & f_{i,j}   \\
  \end{array}\right..
\]
We leave it to the reader to check that $s_{i,j}$ is indeed a homotopy
from $\partial_{i,j}^1\circ\beta_i$ to $\partial_{i,j}^0\circ\beta_j$
and satisfies the cocycle relation \eqref{cocyclerel1}.

It remains to show that $\beta_i$ is a quasi-isomorphism. 
Let $i,j$ be a pair of indices. By definition and by Lemma
\ref{constr2lem} (ii) the restriction of
$\partial_{i,j}^0:\Ll_j\to\Ll_{i,j}$ to $U_{i,j}$ is a quasi-isomorphism.
Therefore, for $q=-1,0$ we may define
\[
   \alpha_{i,j}^q \;:=\; 
     H^q(\partial_{i,j}^0)^{-1}\circ H^q(\partial_{i,j}^1):\; 
           H^q(\Ll_i)|_{U_{i,j}} \iso H^q(\Ll_j)|_{U_{i,j}}.
\]
One checks that the cocycle relation
$\alpha_{j,k}^q\circ\alpha_{i,j}^q=\alpha_{i,k}^q$ holds. Therefore,
there exists a $G$-$\Oo_Y$-module $\T^q$ together with isomorphisms
$\gamma_i^q:\T^q|_{U_i}\iso H^q(\Ll_i)$ such that
$\alpha_{i,j}^q=\gamma_j^q\circ(\gamma_i^q)^{-1}$. It is a bit tedious
but elementary to define, for each $(n+1)$-tuple
$\bi=(i_0,\ldots,i_n)$ an isomorphism
$\gamma_{\bi}^q:\T^q|_{U_{\bi}}\iso H^q(\Ll_{\bi})$ which identifies the
complex
\begin{equation} \label{constr3eq3}
    H^q(\K) \;=\; (\;\prod_i\; H^q(\Ll_i) \;\To\; 
         \prod_{i<j}\; H^q(\Ll_{i,j}) \;\To\; \cdots)
\end{equation}
with the $\check{\rm C}$ech-resolution of the sheaf $\T^q$. We conclude
that the complex \eqref{constr3eq3} is exact.  Now the spectral
sequence $H^p(H^q(\K))\Rightarrow H^{p+q}(\Ll_{Y/S})$ identifies
$\T^q$ with $H^q(\Ll_{Y/S})$ in such a way that $\gamma_i^q$ is
identified with $H^q(\beta_i)$. In
particular, $H^q(\beta_i)$ is an isomorphism, which is what we wanted
to show.  \Endproof

\begin{rem} \label{constrrem1}
  Let $(U_i,\varphi_i)$ and $(U_i',\varphi'_i)$ be two atlases and
  $\Ll_{Y/S}$ and $\Ll_{Y/S}'$ the corresponding complexes, as defined
  above. Then the disjoint union of $(U_i,\varphi_i)$ and
  $(U_i',\varphi_i')$ is again an atlas and gives rise to a third
  complex $\Ll_{Y/S}''$, canonically equipped with quasi-isomorphisms
  $\Ll_{Y/S}''\to\Ll_{Y/S}$ and $\Ll_{Y/S}''\to\Ll_{Y/S}'$. In other
  words: the cotangent complex $\Ll_{Y/S}$, considered as an object of
  the derived category $\Der^+(Y,G)$, does not depend on the choice of
  the atlas $(U_i,\varphi_i)$. 
\end{rem}

\begin{rem} \label{constrrem2}
  By definition we have $H^0(\Ll_{Y/S})=\Omega_{Y/S}$ and
  $H^q(\Ll_{Y/S})=0$ for $q\not\in\{-1,0\}$. It is also clear that
  $\Ll_{Y/S}$ has functorial properties similar to the sheaf of
  differentials $\Omega_{Y/S}$. Namely, if
  \begin{equation} \label{cotangeq2}
  \begin{CD}
    Y'   @>{u}>> Y    \\
    @VVV      @VVV \\
    S'=\Spec R'   @>>> S=\Spec R    \\
  \end{CD}
  \end{equation}
  is a commutative and $G$-equivariant diagram of schemes (where $Y/S$
  and $Y'/S'$ satisfy Assumption \ref{mainass}), then
  we have a natural homomorphism
  \begin{equation} \label{cotangeq3}
      u^*\Ll_{Y/S} \To \Ll_{Y'/S'}
  \end{equation}
  in $\Der^+(Y',G)$. The morphism \eqref{cotangeq3} is an isomorphism in
  each of the following two cases:
\begin{enumerate}
\item
  We have $S=S'$ and $Y'\to Y$ is an open immersion.
\item The diagram \eqref{cotangeq2} is Cartesian and either $Y\to S$
  or $S'\to S$ is flat.
\end{enumerate} 
\end{rem}

\begin{rem} \label{constrrem3}
  If $Y\to S$ is a local complete intersection, then $\Ll_{Y/S}$
  agrees with Illusie's equivariant cotangent complex, up to canonical
  quasi-isomorphism. In general, $\Ll_{Y/S}$ is quasi-isomorphic to
  Illusie's equivariant cotangent complex, truncated at degree $-1$.
  In particular, if $\F$ is a $G$-$\Oo_Y$-module then the $n$th
  hyperext group $\HExt_G^n(\Ll_{Y/S},\F)$ is the `correct' one only
  for $n\leq 1$. See \cite{IllusieCC}, Chapitre III, Corollaire
  1.2.9.1.
\end{rem}


\section{Extensions} \label{ext}

In this section we prove that $G$-equivariant extensions of the
morphism $Y\to S$ by a quasi-coherent $G$-$\Oo_Y$-module $\F$ are
classified by the group $\HExt_G^1(\Ll_{Y/S},\F)$. See Corollary
\ref{extcor}.  This result will be the basis for the results on
equivariant deformations of $Y\to S$ in \S \ref{def}. Actually,
instead of working with extensions of the scheme $Y$, we prefer to
work with the opposite category of extensions of the sheaf $\Oo_Y$.

\subsection{}  \label{ext0}

Let $G\to S=\Spec R$ and $Y\to S$ be as in \S \ref{constr1}. We also
fix a $G$-$\Oo_Y$-module $\F$ which is a quasi-coherent
$\Oo_Y$-module.

\begin{defn}
  An {\em equivariant extension of $\Oo_Y$ by $\F$} is given by a
  short exact sequence of sheaves of $R$-modules on $Y$, of the form
  \[
       0 \;\to\; \F \To \E \To \Oo_Y \;\to\; 0,
  \]
  together with a $G$-action and a structure of sheaf of $R$-algebras
  on $\E$ such that the following holds: 
  \begin{enumerate}
  \item 
    The maps $\F\to\E$ and $\E\to\Oo_Y$ are $G$-equivariant.
  \item 
    The map $\E\to\Oo_Y$ is an $R$-algebra morphism. 
  \item The sheaf $\F$, considered as a subsheaf of $\E$, is
    a sheaf of ideals with square zero. Moreover, the
    induced structure of sheaf of $\Oo_Y$-modules is the
    canonical one.
  \end{enumerate}
  We denote by $\EXT_G(\Oo_Y,\F)$ the category of all such
  extensions. Morphisms between extensions are defined in the obvious
  manner. (The Five Lemma shows that all morphisms are in fact
  isomorphisms.) 
\end{defn}

Given an extension $\E$ of $\Oo_Y$ by $\F$, we get a morphism of
$S$-schemes $Y'\to Y$ (which is a homeomorphism on the underlying
topological spaces) such that $\E=\Oo_{Y'}$. The scheme $Y'$ is called
an extension of $Y$ by $\F$. 

By taking Baer sums of short exact sequences, one defines a bifunctor
$(\E_1,\E_2)\mapsto \E_1+\E_2$. Together with certain natural
transformations $(\E_1+\E_2)+\E_3\cong \E_1+(\E_2+\E_3)$ and
$\E_1+\E_2\cong \E_2+\E_1$, it gives $\EXT_G(\Oo_Y,\F)$ the structure
of a (strictly commutative) Picard category, see Appendix \ref{pic}.

\begin{thm} \label{extthm}
  Let $\F$ be a coherent sheaf of $G$-$\Oo_Y$-modules. We denote by
  $\PIC(\RHom_G(\Ll_{Y/S},\F))$ the Picard category associated to the
  derived complex $\RHom_G(\Ll_{Y/S},\F)$, see Appendix \ref{pic1} and
  \ref{GOY4}. Then there exists a natural isomorphism of Picard
  categories
  \[
     \EXT_G(\Oo_Y,\F) \;\cong\; \PIC(\RHom_G(\Ll_{Y/S},\F)).
  \]
\end{thm}

We will sketch a proof of Theorem \ref{extthm} in the rest
of this section.  The following corollary corresponds to
Theorem 1.5.1 of \cite{IllusieTopos}. It follows from
Theorem \ref{extthm}, using Proposition \ref{pic2prop2}.

\begin{cor} \label{extcor}
  The group of isomorphism classes of equivariant extensions of $Y$ by
  $\F$ is canonically isomorphic to
  $\HExt_G^1(\Ll_{Y/S},\F)$. Moreover, the group of automorphisms of
  any fixed equivariant extension of $Y$ by $\F$ is canonically
  isomorphic to $\Hom_G(\Omega_{Y/S},\F)$.
\end{cor}

\subsection{}  \label{ext1}

In the following three subsections we prove a non-equivariant
version of Theorem \ref{extthm}. To this end, we denote by
$\EXT(\Oo_Y,\F)$ the Picard category of (non-equivariant) extensions of
$\Oo_Y$ by $\F$. 

\begin{prop} \label{ext1prop}
  Suppose that $Y$ is affine and admits a global chart $\varphi:Y\inj
  P$. Let $\Ll_\varphi$ be the cotangent complex of $\varphi$, see \S
  \ref{localconstr}. Then there exists an isomorphism of Picard
  categories
  \[
       F_{\varphi}:\; \EXT(\Oo_Y,\F) \;\liso\; 
         \PIC(\RHom_Y(\Ll_\varphi,\F)).
  \]
  Given a morphism $u:\psi\to\varphi$ of global charts, let
  $\tilde{u}:\PIC(\RHom_Y(\Ll_\psi,\F))\cong
  \PIC(\RHom_Y(\Ll_\varphi,\F))$ denote the isomorphism of Picard
  categories induced from the quasi-isomorphism
  $u^*:\Ll_\varphi\to\Ll_\psi$. There exists an isomorphism of
  additive functors
  \[
       t_u:\; \tilde{u}\circ F_\psi \;\iso\; F_\varphi
  \]
  such that the following holds. If
  $\chi\lpfeil{v}\psi\lpfeil{u}\varphi$ is the composition of two
  morphisms of global charts, then
  \begin{equation} \label{ext1eq1}
        t_{u\circ v} \;=\; t_u\circ\tilde{u}(t_v).
  \end{equation}
\end{prop} 

\proof Using that $Y$ is affine and that $P/S$ is formally smooth one
shows that $\Ll_\varphi^0=\Omega_{P/S}\otimes\Oo_Y$ is a projective
$\Oo_Y$-module. This implies that
\[
     \RHom_Y(\Ll_\varphi,\F)^{[0,1]} \;\cong\; 
         \Hom_Y^\bullet(\Ll_\varphi,\F)^{[0,1]}.
\]
Therefore, we may replace $\RHom_Y(\Ll_\varphi,\F)$ in the statement of
the proposition by the complex
$\Hom_Y^\bullet(\Ll_{Y/S},\F)$. We write $Y=\Spec A$ and
$P=\Spec B$. Then $\varphi$ corresponds to an ideal $I\lhd B$ such
that $A=B/I$. We also write $\F=\widetilde{M}$ for some $A$-module
$M$. With this notation, we have
\[
   \Hom_Y^\bullet(\Ll_\varphi,\F) \;=\; 
     (\; \Hom_A(\Omega_{B/R}\otimes A,M) \;\lpfeil{\circ d}\;
           \Hom_A(I/I^2,M)\;).
\]
An object of $\EXT(\Oo_Y,\F)$ is given by an extension of
$R$-modules $0\to M\to E\to A\to 0$, where $E$ carries in addition the
structure of an $R$-algebra such that the following holds. Firstly,
$E\to A$ is a homomorphism of $R$-algebras; secondly, $M^2=0$,
considered as ideal of $E$. In the rest of the proof, we shall refer
to such an object simply as an {\em extension}. Since $B$ is formally
smooth, there exists a homomorphism of $R$-algebras $\lambda:B\to E$
lifting the canonical map $B\to A$. Set $\nu:=\lambda|_I\mod{I^2}$. It
is clear that $\nu$ is an $A$-linear morphism $I/I^2\to M$. We
consider $\nu$ as an object of $\PIC(\Hom_Y^\bullet(\Ll_\varphi,\F))$
and set
\[
       F_\varphi(E) \;:=\; \nu.
\]

Let $0\to M\to E'\to A\to 0$ be another extension and $f:E\iso E'$ an
isomorphism of extensions. Let $\lambda':B\to E'$ be an $R$-algebra
morphism lifting $B\to A$ and set $\nu':=\lambda'|_I\mod{I^2}$. Then
the map $\lambda'-f\circ\lambda:B\to M$ is easily seen to be an
$R$-linear derivation which vanishes on $I^2$. It corresponds to an
$A$-linear homomorphism $\theta:\Omega_{B/R}\otimes A\to M$ such that
$\theta\circ d=\nu'-\nu$. In other words, $\theta$ is a homomorphism
$\nu\to\nu'$ in $\PIC(\Hom_Y^\bullet(\Ll_\varphi,\F))$. We set
\[
     F_\varphi(f) \;=\; \theta.
\]
One checks that $F_\varphi$ is a faithful additive functor. 

Given an arbitrary $A$-linear map $\nu:I/I^2\to M$, we define the
extension $E_\nu$ as the pushout of $B/I^2$ along $\nu$:
\begin{equation} \label{ext1eq2}
  \begin{array}{ccccccc}
    0 \;\to\;& I/I^2 &\To& B/I^2 &\To& A &\;\to\; 0 \\
    & \rupfeil{\nu} && \rupfeil{\lambda} &&  \rupfeil{=} & \\
    0 \;\to\;& M     &\To& E_\nu    &\To& A &\;\to\; 0 \\
  \end{array}
\end{equation}
It is easy to see that $E_\nu$ carries a unique $R$-algebra structure
such that $M^2=0$ and $\lambda$ is an $R$-algebra morphism. Moreover,
we have $F_\varphi(E_\nu)=\nu$ by construction. Hence $F_\varphi$ is
essentially surjective. 

Let $\theta:\nu\to\nu'$ be an isomorphism in
$\PIC(\Hom_Y^\bullet(\Ll_\varphi,\F))$. This means that $\nu':I/I^2\to
M$ and $\theta:\Omega_{B/R}\otimes A\to M$ are $A$-linear maps such
that $\theta\circ d=\nu'-\nu$. We may identify $\theta$ with the
corresponding derivation $B/I^2\to M$; then
$\theta|_{I/I^2}=\nu'-\nu$. The universal property of the push-forward
shows that there exists a unique $R$-linear map $f:E_\nu\to E_{\nu'}$
such that $f\circ\lambda=\lambda'$. By construction we have
$F_\varphi(f)=\theta$. Hence $F_\varphi$ is fully faithful and even an
isomorphism of Picard categories.

Now let $\psi:Y\inj Q$ be another global chart and $u:\psi\to\varphi$
a morphism of charts. We write $Q=\Spec B'$, $A=B'/{I'}^2$ and
consider $u$ as a morphism of $R$-algebras $B\to B'$. Let $E$ be an
extension and $\lambda:B/I^2\to E$ (resp.\ $\lambda':B'/{I'}^2\to E$)
a lift of $B/I^2\to A$ (resp.\ of $B'/{I'}^2\to A$). By definition we
have
\[
    F_\varphi(E) \;=\; \lambda|_I \pmod{I^2}, \qquad
    \tilde{u}\circ F_\psi(E)\;=\;\lambda'\circ u|_I \pmod{I^2}.
\]
Again it is clear that $\lambda-\lambda'\circ u \pmod{I^2}$ is a
derivation $B/I^2\to M$, corresponding to an $A$-linear map
$\theta:\Omega_{B/R}\otimes A\to M$ and representing a homomorphism
$\tilde{u}\circ F_\psi(E)\to F_\varphi(E)$. We set
\[
       t_u(E) \;=\; \theta.
\]
A formal verification shows that $t_u$ is a morphism of additive
functors $\tilde{u}\circ F_\psi\cong F_\varphi$ and that
\eqref{ext1eq1} holds. 
\Endproof

\subsection{} \label{ext2}

Let $\EEXT(\Oo_Y,\F)$ denote the $Y$-stack whose fiber over a given open
subset $U\subset Y$ is the Picard category $\EXT(\Oo_U,\F|_U)$ (here
$Y$-stack means a stack over the Zariski site of $Y$). It is clear
that $\EEXT(\Oo_Y,\F)$ is a Picard stack, see Appendix
\ref{pic2}.

\begin{prop} \label{ext2prop}
  We assume that $Y$ is affine and admits a global chart
  $\varphi:Y\inj P$.
  \begin{enumerate}
  \item
    Let $U\subset Y$ be an affine open. Then the natural functor
    \begin{equation} \label{ext2eq1}
        \PIC(\RHom_U(\Ll_\varphi|_U,\F|_U)) \;\To\;
           \PPIC(\RHOm_Y(\Ll_\varphi,\F))(U)
    \end{equation}
    is an isomorphism. 
  \item
    There exists a unique isomorphism of Picard stacks
    \[
       \Fu_\varphi:\; \EEXT(\Oo_Y,\F) \;\liso\; 
           \PPIC(\RHOm_Y(\Ll_\varphi,\F))
    \]
    such that for each affine open $U\subset Y$ the restriction of
    $\Fu_\varphi$ to the fiber $\EEXT(\Oo_Y,\F)(U)=\EXT(U,\F|_U)$ is equal
    (up to canonical isomorphism) to the composition of
    $F_{\varphi|_U}$ with \eqref{ext2eq1}.
  \end{enumerate}
\end{prop}

\proof Part (i) follows from Proposition \ref{pic2prop2} and the fact
that the cohomology of the complex $\RHOm_Y(\Ll_\varphi,\F)$ consists
of quasi-coherent sheaves. Part (ii) is a formal consequence of (i)
and is left to the reader.  \Endproof

\subsection{} \label{ext3}

We will now globalize the isomorphism of Picard stacks constructed in
the previous two subsections. To this end, we will use the notation
introduced in \S \ref{globalconstr}. In particular, $(U_i)_{i\in I}$
is a covering of $Y$ by affine opens, admitting local charts
$\varphi_i:U_i\inj P_i$. By Proposition \ref{constr2prop} we
obtain, for each ordered pair $i<j$, an essentially commutative square
of quasi-isomorphisms
\begin{equation} \label{ext3eq1}
\begin{CD}
   \Ll_{Y/S}|_{U_{i,j}}  @>{\beta_i}>> \Ll_{\varphi_i}|_{U_{i,j}} \\
   @V{\beta_j}VV                        @VV{\partial_{i,j}^1}V     \\
   \Ll_{\varphi_j}|_{U_{i,j}} @>>{\partial_{i,j}^0}> \Ll_{\varphi_{i,j}},
\end{CD}
\end{equation}
i.e.\ a homotopy 
$s_{i,j}:\partial_{i,j}^1\circ\beta_i\sim\partial_{i,j}^0\circ\beta_j$
such that the cocycle relation \eqref{cocyclerel1} holds. 
Set
\[
    \PPic \;:=\; \PPIC(\RHOm_Y(\Ll_{Y/S},\F)), \qquad
    \PPic_{\bi} \;=\; \PPIC(\RHOm_Y(\Ll_{\varphi_{\bi}},\F)).
\] 

The diagram \eqref{ext3eq1} yields an essentially commutative square
of isomorphisms of Picard stacks
\begin{equation} \label{ext3eq2}
\begin{CD}
  \PPic_{i,j}   @>{\tilde{\partial}_{i,j}^1}>>  \PPic_i|_{U_{i,j}} \\
  @V{\tilde{\partial}_{i,j}^0}VV   @VV{\tilde{\beta}_i}V   \\
  \PPic_j|_{U_{i,j}}   @>>{\tilde{\beta}_j}>  \PPic|_{U_{i,j}},   \\
\end{CD}
\end{equation}
i.e.\ an isomorphism of additive functors 
$\tilde{s}_{i,j}:\tilde{\beta}_i\circ\tilde{\partial}_{i,j}^1
\sim\tilde{\beta}_j\circ\tilde{\partial}_{i,j}^0$ such that 
\begin{equation} \label{ext3eq3}
   \tilde{\partial}_{i,j,k}^0(\tilde{s}_{j,k}) \circ
     \tilde{\partial}_{i,j,k}^2(\tilde{s}_{i,j}) \;=\;
       \tilde{\partial}_{i,j,k}^1(\tilde{s}_{i,k})
\end{equation}
for all triples $i<j<k$.   

\begin{prop} \label{ext3prop}
  There exists an isomorphism of Picard stacks
  \[
     \Fu:\; \EEXT(\Oo_Y,\F) \;\liso\; 
           \PPic=\PPIC(\RHOm_Y(\Ll_{Y/S},\F))
  \]
  and for each index $i$ an isomorphism of additive functors
  $u_i:\Fu|_{U_i}\cong \tilde{\beta}_i\circ\Fu_{\varphi_i}$.
\end{prop}

\proof
Let $i<j$. By Proposition \ref{ext1prop} we obtain two natural
isomorphisms of additive functors
\[
    F_{\varphi_i}|_{U_{i,j}} \;\cong\; 
      \tilde{\partial}_{i,j}^1 \circ F_{\varphi_{i,j}}, \qquad
    F_{\varphi_j}|_{U_{i,j}} \;\cong\; 
      \tilde{\partial}_{i,j}^0 \circ F_{\varphi_{i,j}}.
\]
Using the essentially commutative square \eqref{ext3eq2} they can be
extended to an isomorphism 
\[
  u_{i,j}:\; \tilde{\beta}_j\circ F_{\varphi_j}|_{U_{i,j}} \cong
   \tilde{\beta}_j\circ\tilde{\partial}_{i,j}^0\circ F_{\varphi_{i,j}}
           \cong
   \tilde{\beta}_i\circ\tilde{\partial}_{i,j}^1\circ F_{\varphi_{i,j}}
           \cong \tilde{\beta}_i\circ F_{\varphi_i}|_{U_{i,j}}.
\]
A tedious but elementary verification, using \eqref{ext1eq1} and
\eqref{ext3eq3}, shows that $u_{i,j}$ satisfies the obvious cocycle
relation. The proposition follows. 
\Endproof

Using the canonical isomorphism
\[
     \RHom_Y(\Ll_{Y/S},\F)) \;\cong\; 
        \RR\Gamma(Y,\RHOm_Y(\Ll_{Y/S},\F))
\]
and Proposition \ref{pic2prop2}, we obtain a non-equivariant version
of Theorem \ref{extthm}:

\begin{cor} \label{ext3cor}
  There is a natural isomorphism of Picard categories
  \[
      F:\; \EXT(\Oo_Y,\F) \;\liso\; \PIC(\RHom_Y(\Ll_{Y/S},\F)).
  \]
\end{cor}

\subsection{} \label{ext4}

We are now going to prove Theorem \ref{extthm} in full generality. In
the sequel, $R'$ will always denote a {\em flat} $R$-algebra, and a
prime stands for base change with respect to $R\to R'$; for instance
$Y':=Y\otimes_R R'$. By Remark \ref{constrrem2} we have natural
isomorphisms
\begin{equation} \label{ext4eq1}
    \Ll_{Y/S}\otimes_R R' \;\cong\; \Ll_{Y'/S'}.
\end{equation}
and
\begin{equation} \label{ext4eq2}
     \RHom_Y(\Ll_{Y/S},\F)\otimes_R R' \;\cong\; 
       \RHom_{Y'}(\Ll_{Y'/S'},\F').
\end{equation}
To simplify the notation, we will henceforth write
\[
     A \;:=\; \RHom_Y(\Ll_{Y/S},\F)^{[0,1]}.
\]
Note that $A$ is a complex of $G$-$R$-modules of amplitude $[0,1]$,
well defined up to canonical isomorphism in $\D^{[0,1]}(Y,G)$.

Let $R'$ be a flat $R$-algebra and $\sigma\in G(R')$. The
automorphism $A'\iso A'$, $a\mapsto a^\sigma$ induces an
isomorphism of Picard categories
$\tilde{\sigma}:\PIC(A')\iso\PIC(A')$. Given an extension
$0\to\F'\to\E'\to\Oo_{Y'}\to 0$ (i.e.\ an object of
$\EXT(\Oo_{Y'},\F)$), let $\E^\sigma$ be the extension
\[
     0 \;\to\; \F'\cong\sigma^*\F \;\To\; \sigma^*\E \;\To\;
          \sigma^*\Oo_{Y'}\cong\Oo_{Y'} \;\to\; 0
\]
Here the isomorphisms $\F'\cong\sigma^*\F$ and
$\sigma^*\Oo_{Y'}\cong\Oo_{Y'}$ come from the $G$-action on $\F$ and
$\Oo_Y$. Given an isomorphism $f:\E_1\iso\E_2$ of extensions, then
$f^\sigma$, as defined in \S \ref{GOY2}, is an isomorphism
$\E_1^\sigma\iso\E_2^\sigma$. One checks that the association
$\E\mapsto\E^\sigma$ is an automorphism of Picard categories
\[
    \tilde{\sigma}:\; \EXT(\Oo_{Y'},\F') \;\liso\; \EXT(\Oo_{Y'},\F').
\]
One checks:

\begin{lem} \label{ext4lem}
  We have an essentially commutative diagram
  \[\begin{CD}
      \PIC(A')     @>{F'}>>    \EXT(\Oo_{Y'},\F)          \\
      @V{\tilde{\sigma}}VV  @VV{\tilde{\sigma}}V \\
      \PIC(A')     @>{F'}>>    \EXT(\Oo_{Y'},\F)          \\
  \end{CD}\]
  (the isomorphism $F'$ is given by Corollary \ref{ext3cor}).
\end{lem}

It follows from Proposition \ref{GRmod3prop} that 
\begin{equation} \label{ext4eq3}
       \RHom_G(\Ll_{Y/S},\F)^{[0,1]} \;\cong\; \Tot(K)^{[0,1]},
\end{equation}
where $K$ is the double complex
\[
     K:\quad\left\{\quad
  \begin{array}{ccc}
    A^0  & \lpfeil{d} & A^1                \\
    \rupfeil{\partial}  & &  \rupfeil{\partial}   \\
    C^1(G,A^0) & \lpfeil{-d} & C^1(G,A^1)  \\
    \rupfeil{\partial}  & &  \rupfeil{\partial}   \\
    C^2(G,A^0) & \lpfeil{d} & C^2(G,A^1)   \\
    \rupfeil{\partial}  & &  \rupfeil{\partial}   \\
      \cdots & & \cdots 
  \end{array}
    \right.
\]
We are now going to construct an isomorphism of Picard categories
\begin{equation} \label{ext4eq4}
  F^G:\; \EXT_G(\Oo_Y,\F) \;\liso\; \PIC(\Tot(K)).
\end{equation}
Together with \eqref{ext4eq3}, this will complete the proof of Theorem
\ref{extthm}.

An object of $\EXT_G(\Oo_Y,\F)$ is an object $\E$ of
$\EXT(\Oo_Y,\F)$, together with an action of $G$ on $\E$ such
that the maps $\F\to\E$ and $\E\to\Oo_Y$ are $G$-equivariant. Such an
action is determined by the following data. For each flat $R$-algebra
$R'$ and group element $\sigma$ we get an isomorphism
$f_\sigma:\E'\iso\E^\sigma$ in $\EXT(\Oo_Y,\F)$ such that
\begin{equation} \label{ext4eq5}
   f_{\sigma\tau} \;=\; f_\sigma^\tau\circ f_\tau
\end{equation}
holds for all pairs $\sigma,\tau\in G(R')$.
Let $F$ be the isomorphism of Corollary \ref{ext3cor} and set
$\nu:=F(\E)\in A^1$, $\theta_\sigma:=F(f_\sigma)\in (A^0)'$. 
By Lemma \ref{ext4lem}, $\theta_\sigma$ is an isomorphism
$\nu'\iso\nu^\sigma$, i.e.\ 
\begin{equation} \label{ext4eq6}
   d(\theta_\sigma) \;=\; \nu^\sigma - \nu'.
\end{equation} 
Equation \eqref{ext4eq5} shows that
\[
    \theta_{\sigma\tau} \;=\; \theta_\sigma^\tau + \theta_\tau.
\]
In other words, the association $\sigma\mapsto\theta_\sigma$
corresponds to a $1$-cocycle $\theta$, i.e.\ an element of
$Z^1(G,A^0)=\Ker(C^1(G,A^0)\pfeil{\partial}C^2(G,A^0))$, see \S
\ref{GRmod3}. Also, Equation \eqref{ext4eq6} means that
$d(\theta)=\partial(\nu)$. We have shown that the pair $(\nu,\theta)$
lies in $Z^1(\Tot(K))$, i.e.\ represents an object of $\PIC(\Tot(K))$.
We set
\[
   F^G(\E) \;:=\; (\nu,\theta).
\]

Now let $g:\E_1\iso\E_2$ be an isomorphism in $\EXT_G(Y,\F)$. Set
$F^G(\E_1):=(\nu_1,\theta_1)$, $F^G(\E_2):=(\nu_2,\theta_2)$ and
$\eta:=F(g)$. Then $d(\eta)=\nu_2-\nu_1$. By definition, $g$ is
$G$-equivariant, so the diagram
\[\begin{CD}
   \E_1'       @>{g'}>>       \E_2'       \\
   @V{f_{1,\sigma}}VV  @VV{f_{2,\sigma}}V \\
   \E_1^\sigma @>{g^\sigma}>> \E_2^\sigma \\
\end{CD}\]
commutes for all $\sigma\in G(R')$. By Lemma \ref{ext4lem} this means
that
\begin{equation} \label{ext4eq7}
  \theta_{2,\sigma} - \theta_{1,\sigma} \;=\; \eta^\sigma - \eta',
\end{equation}
or, equivalently, $\theta_2-\theta_1=\partial(\eta)$. It follows that
$\eta$ corresponds to an isomorphism
$(\nu_1,\theta_1)\iso(\nu_2,\theta_2)$ in $\PIC(\Tot(K))$. We set
\[
    F^G(g) \;:=\; \eta.
\]
We leave it to the reader to check that $F^G$ is indeed an isomorphism
of Picard categories. Now the proof of Theorem
\ref{extthm} is complete.  \Endproof


\section{Deformations}  \label{def}

In this section we show how one can classify equivariant deformations
of $Y\to S$ along an infinitesimal extension $S\inj S'$, using the
equivariant cotangent complex. The main result is Theorem
\ref{defthm}. In \S \ref{def3} we discuss how this result behaves
under localization to a formal neighborhood of a point (Theorem
\ref{def3thm}).

\subsection{} \label{def1}

Let $R'$ be a commutative ring and $\a\lhd R'$ an ideal
with $\a^2=0$.  We set $R:=R'/\a$, $S':=\Spec R'$ and $S:=\Spec R$.
Furthermore, let $G'\to S'$ be a flat affine group scheme and $Y\to S$
a flat morphism together with an $S$-linear action of
$G:=G'\times_{S'} S$ on $Y$. 

\begin{defn} \label{def1def}
  An {\em equivariant deformation} of $Y\to S$ to $S'$ is a flat
  morphism $Y'\to S'$ together with an $S'$-linear action of $G'$ on
  $Y'$ and a $G$-equivariant isomorphism of $S$-schemes $Y\cong
  Y'\times_{S'} S$. An isomorphism of deformations is a
  $G'$-equivariant isomorphism of $S'$-schemes $Y_1'\cong Y_2'$ which
  induces the identity on $Y$. 
\end{defn}

Theorem \ref{defthm} below shows how to classify isomorphism classes
of equivariant deformations of $Y\to S$ using the equivariant
cotangent complex $\Ll_{Y/S}$. However, in the proof of Theorem
\ref{defthm} we will also use the cotangent complex of the composed
morphism $Y\to S\inj S'$. Therefore, to be able to use the definition
of $\Ll_{Y/S'}$ in \S \ref{constr}, we make the following assumption.

\begin{ass} \label{mainass2}
  Every point of $Y$ is contained in an affine and $G$-stable open
  $U\subset Y$ such that the following holds. There exists a smooth
  affine $S'$-scheme $P'\to S'$, an $S'$-linear action of $G'$ on $P'$
  and a $G$-equivariant closed immersion $\varphi:U\inj P'$. In other
  words, Assumption \ref{mainass} holds for the composed morphism
  $Y\to S'$ and the group scheme $G'$.
\end{ass}

Under this assumption we can prove:

\begin{thm} \label{defthm}
  \begin{enumerate}
  \item
    There exists an element (called the obstruction)
    \[
       \omega=\omega(Y/S,S') \;\in\; 
          \HExt_G^2(\Ll_{Y/S},\Oo_Y)\otimes_R\a,
    \]
    depending functorially on $Y\to S\inj S'$, whose vanishing is
    necessary and sufficient for the existence of an equivariant
    deformation of $Y\to S$ to $S'$.
  \item
    Suppose that $\omega=0$. Then the set of isomorphism classes of
    deformations of $Y\to S$ to $S'$ is, in a natural way, a principal
    homogeneous space under the abelian group
    \[
         \HExt_G^1(\Ll_{Y/S},\Oo_Y)\otimes_R\a.
    \]
  \end{enumerate}
\end{thm}

This corresponds to Proposition 2.3 of \cite{IllusieTopos}. However,
if $Y\to S$ is not a local complete intersection, then our definition
of $\Ll_{Y/S}$ does not always give the same hyperext group
$\HExt_G^2(\Ll_{Y/S},\Oo_Y)$ as in \cite{IllusieCC} and
\cite{IllusieTopos}. In particular, our obstruction does not
necessarily agree with the obstruction constructed in
\cite{IllusieCC} and \cite{IllusieTopos}, simply because it does not
lie in the same group. See also Remark \ref{def2rem} below.

\subsection{Proof of Theorem \ref{defthm}} \label{def2}

Let $\F:=\Oo_Y\otimes_R \a$. Since $Y\to S$ is flat we have natural
isomorphisms 
\begin{equation} \label{def1eq1}
   \HExt_G^n(\Ll_{Y/S},\F) \;\cong\; 
       \HExt_G^n(\Ll_{Y/S},\Oo_Y)\otimes_R\a
\end{equation}
for all $n$. By Assumption \ref{mainass2} and Definition \ref{Ldef},
the equivariant cotangent complexes $\Ll_{Y/S}$ and $\Ll_{Y/S'}$
are defined as complexes of $G$-$\Oo_Y$-modules and we have a natural
$G$-equivariant morphism $\Ll_{Y/S'}\to\Ll_{Y/S}$.

\begin{lem} \label{def1lem}
  There is a natural exact sequence in $\Der^+(Y,G)$
  \begin{equation} \label{def1eq2}
      0 \;\To\; \F[1] \;\To\; \Ll_{Y/S'} \;\To\; \Ll_{Y/S} 
        \;\To\; 0.
  \end{equation}
  More precisely, the natural morphism $\Ll_{Y/S'}\to \Ll_{Y/S}$
  is surjective in all degrees, and there exists a $G$-equivariant
  quasi-isomorphism $\F[1]\to\Ker(\Ll_{Y/S'}\to \Ll_{Y/S})$.
  (Recall that $\F[1]$ denotes the complex where $\F$ is placed in
  degree $-1$.)
\end{lem}

\proof Let $\varphi':U\inj P'$ be a local chart for the morphism $Y\to
S'$ and $\I'\subset\Oo_{P'}$ the corresponding sheaf of ideals. Then
$\varphi'$ gives rise to a local chart $\varphi:U\inj
P:=P\times_{S'}S$ for the morphism $Y\to S$. The corresponding sheaf
of ideals is $\I:=\I'/\F$. It is clear that
$\Omega_{P'/S'}\otimes\Oo_Y\cong\Omega_{P/S}\otimes\Oo_Y$. Moreover, we
have a short exact sequence
\begin{equation} \label{def1eq3}
  0 \;\To\; \F \;\To\; \I'/{\I'}^2 \;\To\; \I/\I^2 \;\To\; 0.
\end{equation}
Hence it follows from Definition \ref{Ldef} that
$\Ll_{Y/S'}\to\Ll_{Y/S}$ is surjective in all degrees and that its
kernel is isomorphic to the $\check{\rm C}$ech-resolution of $\F[1]$
(with respect to the open covering $(U_i)$ used to define
$\Ll_{Y/S}$). This proves the lemma.
\Endproof

Let $\G$ be a $G$-$\Oo_Y$-module.  The short exact sequence
\eqref{def1eq2} gives rise to the following long exact sequence
\begin{equation} \label{def1eq4}
  0\;\to\; \HExt_G^1(\Ll_{Y/S},\G) \;\To\; \HExt_G^1(\Ll_{Y/S'},\G)
    \;\To\; \Hom_G(\F,\G) \;\lpfeil{\partial} \HExt_G^2(\Ll_{Y/S},\G).
\end{equation}
This applies in particular to the case $\G:=\F$.  We define the
obstruction $\omega:=\omega(Y/S,S')$ as the image of the identity map
$\Id:\F\to\F$ under the boundary map $\partial$ in
\eqref{def1eq4}. Now Theorem \ref{defthm} follows from Corollary
\ref{extcor}, the exactness of \eqref{def1eq4} and the following
proposition.

\begin{prop} \label{def1prop}
  There is a natural bijection between 
  \begin{itemize}
  \item[(a)]
    deformations of $Y\to S$ to $S'$, up to isomorphism, and
  \item[(b)] elements of $\HExt_G^1(\Ll_{Y/S'},\F)$ which are mapped to
    $\Id_{\F}$ (by the middle arrow in \eqref{def1eq4}).
  \end{itemize}
\end{prop}

\proof 
Let $\G$ be a quasi-coherent $G$-$\Oo_Y$-module. 
By Corollary \ref{extcor}, an element of
$\HExt_G^1(\Ll_{Y/S'},\G)$ corresponds to an equivariant extensions of
$Y$ by $\G$, i.e.\ a closed equivariant embedding $Y\inj Y'$ of
$S'$-schemes defined by an ideal $\J\subset\Oo_{Y'}$, together with an
isomorphism $\G\cong\J$ of $G$-$\Oo_Y$-modules. We obtain a morphism of
$G$-$\Oo_Y$-modules
\begin{equation} \label{def1eq5}
     \F=\Oo_Y\otimes_R\a \;\To\; \J\cong\G.
\end{equation}
By reexamination of the proof of Theorem \ref{extthm} one shows that
the middle arrow of the sequence \eqref{def1eq4} maps the element of
$\HExt_G^1(\Ll_{Y/S'},\G)$ corresponding to the extension $Y'$ to the
morphism \eqref{def1eq5}.  Also, the local criterion of flatness (see
\cite{Matsumura}, Theorem 49) shows that the morphism $Y'\to S'$ is
flat if and only if \eqref{def1eq5} is an isomorphism. 

The proposition follows easily from these arguments. First, an
equivariant extension of $Y$ by $\F$ for which \eqref{def1eq5} is the
identity on $\F$ gives rise to an equivariant deformation of $Y\to S$
to $S'$. Conversely, let $Y'\to S'$ be an equivariant deformation of
$Y\to S$, and let $\J\subset\Oo_{Y'}$ be the sheaf of ideals
corresponding to the embedding $Y\inj Y'$. Since $Y'\to S'$ is flat by
assumption, the natural map $\F\to \J$ is an isomorphism. Using this
isomorphism, we can see $Y'$ as an equivariant extension of $Y$ by
$\F$ for which \eqref{def1eq5} is the identity on $\F$.  This
concludes the proof of the proposition and hence of Theorem
\ref{defthm}.  \Endproof

\begin{rem}  \label{def2rem}
  The short exact sequence of Lemma \ref{def1lem} should be compared
  with the {\em transitivity triangle} attached to the composition of
  morphisms 
  $Y\to S\inj S'$ in \cite{IllusieCC}:
  \begin{equation} \label{def1eq6}
   \begin{array}{ccccc}
     & & \Ll_{Y/S}^I  & &  \\
      & \!\!\swarrow\!\!\!\! & & \!\!\!\!\nwarrow\!\! & \\
     \Ll_{S/S'}^I\otimes\Oo_Y\!\!\!\! & & \To & & \Ll_{Y/S'}^I \\
  \end{array}
  \end{equation}
  Here $\Ll_{Y/S}^I$ denotes the cotangent complex in the sense of
  Illusie. 
  
  Now suppose that $Y\to S$ is a local complete intersection. Then
  $\Ll_{Y/S}\cong\Ll_{Y/S}^I$. We also have natural morphisms
  $\Ll_{Y/S'}^I\to\Ll_{Y/S'}$ and $\Ll_{S/S'}^I\otimes\Oo_Y\to\F[1]$, but
  they are quasi-isomorphisms only if $S\inj S'$ is a local complete
  intersection (which is typically not the case). Nevertheless, one
  can show that the obstruction $\omega$ in Theorem \ref{defthm} is
  the same as the obstruction obtained by Illusie's theory (via the
  canonical isomorphism
  $\HExt_G^2(\Ll_{Y/S},\Oo_Y)\cong\HExt_G^2(\Ll_{Y/S}^I,\Oo_Y)$).
\end{rem}

\subsection{Localization} \label{def3}

Keeping the notation introduced before, we now impose the following
finiteness conditions.

\begin{ass} \label{def3ass}
\begin{enumerate}
\item
  The affine scheme $S=\Spec R$ is local, Artinian and Noetherian.
\item
  The group scheme $G$ is finite and flat over $S$.
\item
  The scheme $Y$ is either of finite type over $S$ or the localization
  of something of finite type over $S$.
\end{enumerate}
\end{ass}
It follows from Part (i) and (iii) of the assumption that $Y$ is
Noetherian.

By Assumption \ref{mainass} the action of $G$ on $Y$ is admissible;
hence the quotient scheme $X:=Y/G$ exists. It follows from Assumption
\ref{def3ass} (ii) that the projection $\pi:Y\to X$ is finite. Let
$x\in X$ be a point. Let $\Xd=\Spec \Od_{X,x}$ denote the completion
of $X$ at $x$ and set $\Yd:=Y\times_X\Xd$. Since $\pi:Y\to X$ is
finite, $\Yd$ is naturally isomorphic to the completion of $Y$ along
the fiber $\pi^{-1}(x)$. The action of $G$ on $Y$ induces an action of
$G$ on $\Yd$. Since $\Xd\to X$ is flat, we have $\Yd/G=\Xd$.

Let $u:\Yd\to Y$ denote the canonical map. By Remark \ref{constrrem2}
we have a canonical morphism of complexes of $G$-$\Oo_{\Yd}$-modules
\begin{equation} \label{def3eq1}
  u^*\Ll_{Y/S} \;\To\; \Ll_{\Yd/S}.
\end{equation}
A technical complication arises from the fact that \eqref{def3eq1} is
in general not a quasi-isomorphism. However, the next proposition
shows that this does not really matter to us. 

\begin{prop} \label{def3prop}
  Let $\F$ be a coherent sheaf of $G$-$\Oo_Y$-modules.
  There exists an isomorphism of Picard categories
  \[ 
       F_x:\,\EXT_G(\Oo_{\Yd},u^*\F) \;\liso\;
          \PIC(\RHom_G(u^*\Ll_{Y/S},u^*\F))
  \]
  such that the following diagram commutes:
  \begin{equation} \label{def3eq3}
  \begin{CD}
     \EXT_G(\Oo_Y,\F) @>>> \EXT_G(\Oo_{\Yd},u^*\F) \\
     @V{F}VV                    @VV{F_x}V                  \\
     \PIC(\RHom_G(\Ll_{Y/S},\F))  @>>>      
                       \PIC(\RHom_G(u^*\Ll_{Y/S},u^*\F)).    \\
  \end{CD}
  \end{equation}
  Here the upper horizontal arrow is the functor which sends an
  extension $Y'$ of $Y$ by $\F$ to the completion of $Y'$ along the
  fiber $\pi^{-1}(x)$. The left vertical arrow is the isomorphism from
  Theorem \ref{extthm}. The lower horizontal arrow is the natural
  pullback map. 
\end{prop}

\proof The construction of the equivalence $F_x$ is very similar to
the construction of $F$ in the proof of Theorem \ref{extthm}. An
essential difference appears only in the first step, see \S
\ref{ext1}. We will therefore assume for the rest of the proof that
$Y$ is affine and that $G=1$.

Replacing $Y$ by its localization at any point $y\in\pi^{-1}(x)$ we
may assume that $Y=\Spec A$ is local. We have $\Yd=\Spec\Ad$, where
$\Ad$ is the completion of $A$.  The coherent sheaf $\F$ is given by a
finite $A$-module $M$; the pullback $u^*\F$ corresponds to the
$\Ad$-module $\Mh:=M\otimes_A\Ad$. Since $M$ is a finite $A$-module,
$\Mh$ is the $\m_A$-adic completion of $M$.  By Assumption
\ref{def3ass} (iii) we can write $A=B/I$, where $B$ is the
localization of a polynomial ring over $R$ and $I\lhd B$ is an ideal.
By Assumption \ref{def3ass} (i) the ring $B$ is Noetherian and hence
$I$ is finitely generated. Moreover, $\Ad=\Bd/\Ih$, where $\Bd$ is the
completion of $B$ at its maximal ideal and $\Ih:=I\Bd$. The ring $\Bd$
is a power series ring over $R$. In general, $B$ is not formally
smooth over $R$ but only $\m_B$-smooth (see \cite{MatsumuraII}; note
that `formal smoothness' is called `$0$-smoothness' in loc.cit.).

The complex $u^*\Ll_{Y/S}$ corresponds to the
complex of $\Ad$-modules
\[
   \Ld \;:=\; (\;\Ih/\Ih^2 \;\To\; \Omega_{B/R}\otimes_B \Ad\;).
\]
The canonical map
$\Omega_{B/R}\otimes_B\Ad\to\Omega_{\Bd/R}\otimes_{\Bd}\At$ is
injective but in general not surjective. However,
$\Omega_{B/R}\otimes\Ad$ is mapped isomorphically onto
$\Omega_{\Bd/R}^{\rm cont}\otimes\Ad$, where 
\[
    \Omega_{\Bd/R}^{\rm cont} \;:=\; 
      \Omega_{\Bd/R}/(\cap_n\, \m_{\Bd}^n\cdot\Omega_{\Bd/R})
\]
denotes the module of {\em continuous differentials}.  Since
$\Omega_{B/R}\otimes\Ad$ is a free $\Ad$-module, we have
\begin{equation} \label{def3eq4}
    \RHom_Y(u^*\Ll_{Y/S},u^*\F)^{[0,1]} \:=\;
       (\,\Hom_{\Ad}(\Omega_{B/R}\otimes\Ad,\Mh) \;\To\;
          \Hom_{\Ad}(\Ih/\Ih^2,\Mh)\,).
\end{equation}
An object of the Picard category $\EXT(\Oo_{\Yd},u^*\F)$ is given by an
extension $\Mh\to E\to\Ad$ of $R$-modules, with $\Mh^2=0$. In what
follows we will refer to such an object simply as an {\em extension}. 

\begin{lem}
  Let $\Mh\to E\to\Ad$ be an extension and let $\m_E\lhd E$
  denote the inverse image of the maximal ideal $\m_{\Ad}$ of $\Ad$.
  \begin{enumerate}
  \item
    The ring $E$ is complete with respect to the ideal $\m_E$.
  \item
    There exists a continuous lift $\lambda:\Bd\to E$ of the canonical
    map $\Bd\to\Ad$. 
  \end{enumerate}
\end{lem}

\proof
Look at the following ladder with exact rows:
\begin{equation} \label{def3eq5}
\begin{array}{ccccccccc}
  0 &\;\to\;& \Mh &\;\To\;& E &\;\To\;& \Ad &\;\to\;& 0  \\
    &&  \downarrow  &&  \downarrow  && \downarrow  &&    \\
  0 &\;\to\;& \varprojlim \Mh/(\m_E^n\cap\Mh) &\;\to\;&
              \varprojlim E/\m_E^n &\;\to\;& 
              \varprojlim \Ad/\hat{\m}^n  &\;\to\;&  0.
\end{array}
\end{equation}
The vertical arrow on the right is an isomorphism by definition.  An
argument similar to the one used in the proof of the Artin--Rees Lemma
(see \cite{MatsumuraII}, Theorem 8.5) shows that there exists a
constant $c>0$ such that
\[
       \m_E^n\cap\Mh \;\subset\; \m_{\Ad}^{n-c}\cdot\Mh
\]
for all $n\geq c$ (here we use that $A$ and hence $\Ad$ is
Noetherian). Therefore, the vertical arrow in \eqref{def3eq5} on the
left is an isomorphism. Now the Five-Lemma implies that the vertical
arrow in the middle is an isomorphism, too. This proves (i). Part (ii)
of the lemma follows from Part (i) and the $\m_B$-smoothness of $\Bd$.
\Endproof

Using this lemma, the construction of the equivalence $F_x$ is
essentially the same as in the proof of Proposition \ref{ext2prop}. It
is also clear from this construction that the diagram \eqref{def3eq3}
commutes. There are two points one has to pay attention to. The first
is to consider only continuous lifts $\lambda:\Bd\to E$. The second is
this: if $E'$ is another extension, $\lambda':\Bd\to E'$ a lift and
$f:E\iso E'$ an isomorphism of extensions, then
$\lambda'-f\circ\lambda:\Bd\to\Mh$ is a {\em continuous} $R$-linear
derivation which vanishes on $\Ih^2$; it therefore corresponds to an
$\Ad$-linear map $\theta:\Omega_{B/R}\otimes\Ad\to\Mh$. Moreover, any
$R$-linear derivation $\Bd\to\Mh$ is automatically continuous because
$\Mh$ is complete and hence separated with respect to the $\m_A$-adic
topology. This completes the proof of the proposition.  \Endproof

\begin{rem} 
  The proposition is essentially equivalent with the statement that
  the homomorphism
  \[
      \HExt_G^n(\Ll_{\Yd/S},u^*\F) \;\To\; 
         \HExt_G^n(u^*\Ll_{\Yd/S},u^*\F)
  \]
  induced by \eqref{def3eq1} is an isomorphism for $n=0,1$. I suspect
  that this is true for $n>1$ as well, but I don't know how to prove
  this. 
\end{rem}

For $n\geq 0$ we write $\EXt_G^n(\Ll_{Y/S},\F)\sphat_x$ for the
$\Od_{X,x}$-module $\EXt_G^n(\Ll_{Y/S},\F)_x\otimes\Od_{X,x}$. It
follows from flatness of $\Xd\to X$ that
\[
   \EXt_G^n(\Ll_{Y/S},\F)\sphat_x \;=\;
     \HExt_G^n(u^*\Ll_{Y/S},u^*\F).
\]
The local-global spectral sequence from \S \ref{GOY5} gives rise to a
  {\em localization map}
\[
     \HExt_G^n(\Ll_{Y/S},\F) \;\To\;
         \EXt_G^n(\Ll_{Y/S},\F)\sphat_x.
\]

\begin{thm} \label{def3thm}
  Let $S\inj S'=\Spec R'$ be a small extension, with $R=R'/\a$. Let
  $\omega$ be the obstruction for
  lifting $Y$ to $S'$. Also, let $\omega_x$ denote the image of
  $\omega$ under the localization map 
  \[
     \HExt_G^2(\Ll_{Y/S},\Oo_Y)\otimes\a \;\To\;
         \EXt_G^2(\Ll_{Y/S},\Oo_Y)\sphat_x\otimes\a.
  \]
  \begin{enumerate}
  \item
     There exists an equivariant deformation of $\Yd$ to $S'$ if and
     only if $\omega_x=0$.
  \item
     If $\omega_x=0$ then the set of isomorphism classes of
     deformations of $\Yd$ to $S'$ is a principal homogeneous space
     under the group $\EXt_G^1(\Ll_{Y/S},\Oo_Y)\sphat_x\otimes\a$.
  \item
     If $\omega=0$ then the action of $\HExt_G^1(\Ll_{Y/S},\F)$ on the
     set of isomorphism classes of deformations of $Y$ to $S'$ is
     compatible with the action of
     $\EXt_G^1(\Ll_{Y/S},\Oo_Y)\sphat_x\otimes\a$ on deformations of
     $\Yd$, with respect to the localization map
     \[
         \HExt_G^1(\Ll_{Y/S},\Oo_Y)\otimes\a \;\To\;
              \EXt_G^1(\Ll_{Y/S},\Oo_Y)\sphat_x\otimes\a.
     \]
  \end{enumerate}
\end{thm}

\proof
This is proved in the same way as Theorem \ref{defthm}, except that
the exact sequence \eqref{def1eq4} is replaced by the sequence
\begin{equation} \label{def3eq8}
  0\,\to\, \EXt_G^1(\Ll_{Y/S},\F)\sphat_x \;\to\; 
     \EXt_G^1(\Ll_{Y/S'},\F)\sphat_x
    \;\to\; \HOm_G(\F,\F)\sphat_x \;\pfeil{\partial} 
       \EXt_G^2(\Ll_{Y/S},\F)\sphat_x
\end{equation}
and we use Proposition \ref{def3prop} in addition to Theorem
\ref{extthm}. The compatibility statement (iii) follows from the
commutativity of the diagram \eqref{def3eq3} and the fact
that the localization maps define a homomorphism between the exact
sequences \eqref{def1eq4} and \eqref{def3eq8}.  \Endproof


\section{Multiplicative deformation data}  \label{defdat}

Let $X$ be a smooth projective curve, defined over an algebraically
closed field $k$ of characteristic $p>0$.  A multiplicative
deformation datum over $X$ is a pair $(Z,V)$, where $Z\to X$ is a
Galois cover, with Galois group $H$ of order prime to $p$, and an
$H$-stable $\FF_p$-vector space $V$ of logarithmic differential forms
on $Z$. In \S \ref{defdat1}, we associate to the pair $(Z,V)$ a
singular curve $Y$ together with an action of a finite group scheme
$G$ such that $X=Y/G$.  Essentially, $G$ is a semi-direct product
$\bmu_p^s\rtimes H$ (where $s:=\dim_{\FF_p}V$) and $Y\to Z$ is
generically a $\bmu_p^s$-torsor determined by a basis
$\phi_1,\ldots,\phi_s$ of $V$.

As an application of the general theory developed in the previous
sections, we study equivariant deformations of $Y$. Even though
the cover $Y\to X$ is inseparable, its deformation theory is in some
sense similar to the deformation theory of a tame cover. For instance,
we get a morphism of deformation functors
\[
          \Def(Y,G) \;\To\; \Def(X;\tau_j),
\]
see \S \ref{defdat2} for a precise definition. In \S \ref{defdat3} we
give a criterion when this morphism is an isomorphism. The reason for
this relatively nice behavior of $\Def(Y,G)$ is that the `$p$-Sylow
subgroup' of $G$ is a {\em multiplicative} group scheme, whose
cohomology is trivial.  Thus, all the contribution to the hyperext
groups $\HExt_G^n(\Ll_{Y/k},\Oo_Y)$ comes from the cohomology of a
certain coherent sheaf on $X$, and there is no group cohomology
involved.

\subsection{The $G$-cover associated to a deformation
  datum} \label{defdat1}

Fix an algebraically closed field $k$ of characteristic $p>0$ and a
smooth $k$-curve $X$. Let $H$ be a finite group of prime-to-$p$ order
and $\chi$ a character of $H$ with values in $\FF_p$.

\begin{defn} \label{defdatdefn}
  A {\em (multiplicative) deformation datum} on $X$ of type $(H,\chi)$
  is a pair $(Z,V)$, where
  \begin{itemize}
  \item $\pi:Z\to X$ is a finite, tamely ramified Galois cover with
    Galois group $H$, and
  \item
    $V\subset \Omega_{k(Z)/k}$ is an $H$-stable and finite dimensional
    $\FF_p$-vector space consisting of logarithmic differential
    forms on $Z$.  Let $V_k$ denote the $k$-linear hull of $V$ in
    $\Omega_{k(Z)/k}$.  We demand that $\dim_k V_k=\dim_{\FF_p}V$ and
    that $H$ acts on $V$ with character $\chi$. 
  \end{itemize} 
  Recall that a differential form $\phi\in\Omega_{k(Z)/k}$ is called
  {\em logarithmic} if it can be written as $\phi=\diff u/u$ for some
  rational function $u\in k(Z)$.  
\end{defn}

If $\dim_{\FF_p}V=1$ then Definition \ref{defdatdefn} agrees with
Definition 1.5 of \cite{bad}. In this paper we shall only consider
multiplicative deformation data (as opposed to {\em additive}
deformation data), so we omit from now on the adjective
`multiplicative'.

Let us fix a deformation datum $(Z,V)$ of type $(H,\chi)$. For the
moment, we will consider $V$ simply as a (right) $\FF_p[H]$-module.
Let $W(k)$ denote the ring of Witt vectors over $k$ and $W(k)[V]$ the
group ring of $V$ over $W(k)$ (here we consider $V$ as an abelian
group).  Then
\[
        G_0 \;:=\; \Spec W(k)[V]
\]
is a finite flat and commutative group scheme over $W(k)$. In fact, $G_0$
represents the group functor (on the category of $W(k)$-algebras)
\[
    R \;\longmapsto\; G_0(R)=\Hom_{\rm gr}(V,R^\times).
\]
Groups schemes of this form are called {\em diagonalizable} in
\cite{SGA3}, Expos\'e I. We shall write
$\underline{\zeta}=(\zeta_\phi)_{\phi\in V}$ for an element of
$G_0(R)$. Here $\zeta_\phi\in R^\times$ such that
$\zeta_{\phi_1}\zeta_{\phi_2}=\zeta_{\phi_1+\phi_2}$. In particular,
$\zeta_\phi^p=1$. Therefore, the choice of an $\FF_p$-basis of $V$
gives rise to an isomorphism $G_0\cong \bmu_p^n$, where
$n=\dim_{\FF_p}V$.

An element $\beta\in H$ induces an automorphism $G_0\iso G_0$ of group
schemes which sends $\underline{\zeta}=(\zeta_\phi)_{\phi\in V}\in
G_0(R)$ to
\[
   \beta(\underline{\zeta}) \;:=\; (\phi\mapsto\zeta_{\beta^*\phi})
      \;\in G_0(R).
\]
This gives an action of $H$ on $G_0$ from the left. We define the
group scheme $G$ as the semidirect product $G_0\rtimes H$; it
represents the group functor
\[
   R \;\longmapsto\; G(R):=G_0(R)\rtimes H.
\]
The multiplication on the right hand is determined by the rule
\[
   (\underline{\zeta}_1,\beta_1)\cdot(\underline{\zeta}_2,\beta_2)
     \;:=\; (\underline{\zeta}_1\cdot\beta_1(\underline{\zeta}_2),
                \beta_1\beta_2).
\]
Note that the subgroup scheme $G_0\subset G$ is equal to the local part
of $G$. 

Let $R$ be a $W(k)$-algebra and $M$ a $G$-$R$-module. The induced
action of $G_0$ on $M$ is given by a map $\mu:M\to R[V]\otimes_R
M$. It gives rise to a {\em $V$-grading},
i.e.\ a direct sum decomposition
\[
    M \;=\; \bigoplus_{\phi\in V} \, M_\phi, \qquad
       M_\phi \;:=\; \{\,m\in M\mid \mu(m) = \phi\otimes m\,\}. 
\]
One checks that a $G$-$R$-module is the same as an $R$-module together
with a $V$-grading and an $R$-linear action of $H$ from the right such
that $M_\phi^\beta=M_{\beta^*\phi}$ for all $\beta\in H$ and $\phi\in
V$. See also \cite{SGA3}, Expos\'e I. Using the assumption that the
order of $H$ is prime to $p$ one shows:

\begin{lem} \label{defdat1lem}
  Let $R$ be a $W(k)$-algebra and $M$ a $G$-$R$-module. Then
  \[
      H^n(G,M) \;= \quad
      \begin{cases} 
        \quad M_0^H  & \text{\rm for $n=0$},\\
         \quad 0     & \text{\rm for $n>0$}.
      \end{cases}
  \]
\end{lem}

\begin{constr}  \label{Yconstr}
  Let $(Z,V)$ be a deformation datum of type $(H,\chi)$ over $X$.
  We shall construct a curve $Y$ over $k$ and a $G$-action on $Y$ such
  that $Z=Y/G_0$ and $X=Y/G$. The definition of $Y$ and the $G$-action
  will depend, up to canonical isomorphism, only on the deformation
  datum $(Z,V)$ but not on the choices we make during the
  construction. Therefore, it suffices to give the construction
  locally on $X$. Hence, we may assume that $Z=\Spec A$ is affine. Let
  us also choose an $\FF_p$-basis $\phi_1,\ldots,\phi_s$ of $V$. Since
  $\phi_i$ is logarithmic, we have $\phi_i=\diff u_i/u_i$ for some
  rational function $u_i$ on $Z$. After shrinking $Z$ and replacing
  $\phi_i$ by a suitable $\FF_p$-multiple of itself, we may assume
  that $u_i$ lies in $A$ and has at most simple zeros on $Z$. Set
  \[
     B \;:=\; A[y_1,\ldots,y_s\mid y_i^p=u_i], \qquad
     Y \;:=\; \Spec B.
  \]
  The $A$-algebra $B$ has a unique $V$-grading such that $B_0=A$ and
  $y_i\in B_{\phi_i}$. It gives rise to an action of $G_0$ on $Y$ such
  that $Z=Y/G_0$. One checks that there is a unique way to extend the
  action of $H$ on $A$ to an action on $B$ such that
  $\beta^*B_\phi=B_{\beta^*\phi}$. Whence an action of $G$ on $Y$ such
  that $Z=Y/G_0$ and $X=Y/G$. This finishes the construction of $Y$.
\end{constr}

\begin{defn} \label{Ydefn}
  Let $(Z,V)$ be a deformation datum of type $(H,\chi)$, and let $Y$
  be the $k$-curve with $G$-action from Construction \ref{Yconstr}.
  Let $\tau\in X$ be a closed point and choose a point $\xi\in Z$
  above $\tau$. We say that $\tau$ is
  \begin{enumerate}
  \item
    a {\em tame branch point} if it is a branch point of the tame
    cover $Z\to X$, 
  \item
    a {\em wild branch point} if there exists $\phi\in V$ such that
    $\ord_\xi\phi=-1$,
  \item
    a {\em critical point} if it is a branch point (tame or wild) or
    if 
    \[
         \min_{\phi\in V}\,(\,\ord_\xi\phi\,) \;\not=\; 0.
    \]
  \end{enumerate}
  Note that these conditions do not depend on the choice of $\xi$ and
  that a branch point can be wild and tame at the same time.
\end{defn}

\begin{notation} \label{Ynot}
  Let $(\tau_j)_{j\in B}$ denote the set of critical points for
  $(Z,V)$, indexed by the finite set $B$. Let $B\tame$ (resp.\ 
  $B\wild$) denote the subset of $B$ corresponding to the tame (resp.\ 
  wild) branch points; set $B\branch:=B\tame\cup B\wild$. We have a
  divisor on $Z$
  \[
     D \;:=\; \sum_{\xi\in Z}
          \;(\,\min_{\phi\in V}\,\ord_\xi\phi\,) \cdot \xi.
  \]
  We can write $D$ as the difference of two disjoint effective divisors in a
  unique way:
  \[
       D \;=\; D_0 - D_\infty.
  \]
  Note that the image of $D$ (resp.\ of $D_\infty$) on $X$ has support
  in the set of critical points (resp.\ in the set of wild branch
  points). 
\end{notation}

\begin{rem}
\begin{enumerate}
\item
  The map $Y\to X$ is finite and flat. It is a $G$-torsor precisely
  outside the set of branch points. 
\item The curve $Y$ is generically smooth over $k$ if and only if
  $\dim_{\FF_p}V=1$. If this is the case then the singular
  points of $Y$ are precisely the points lying over a critical point
  which is not a wild branch point. 
\end{enumerate}
\end{rem}

\subsection{Equivariant deformations of $Y$} \label{defdat2}

Let $\C_k$ denote the category of local Artinian $W(k)$-algebras.  A
{\em $G$-equivariant deformation} of $Y$ over $R\in\C_k$ is a flat
$R$-scheme $Y_R$ together with an action of $G$ and a $G$-equivariant
isomorphism $Y\cong Y_R\otimes k$ (compare with Definition
\ref{def1def}). We are concerned with the deformation functor
\[
      R \;\longmapsto\; \Def(Y,G)(R)
\]
which sends $R\in\C_k$ to the set of isomorphism classes of
$G$-equivariant deformations of $Y$ over $R$. 
The next lemma follows easily from Construction \ref{Yconstr}: 

\begin{lem}
  Let $Y_R$ be an equivariant deformation of $Y$ over
  $R$. Furthermore, let $R'\to R$ be a small extension, i.e.\
  $R=R'/\a$ for an ideal $\a\lhd R'$ such that
  $\a\cdot\m_{R'}=0$. Then the morphism $Y_R\to\Spec R'$ satisfies
  Assumption \ref{mainass}.
\end{lem}

The lemma shows that the equivariant cotangent complex $\Ll_{Y/k}$ is
defined and that we can apply Theorem \ref{defthm} to classify the set
of liftings of the deformation $Y_R$ to $R'$. Let $k[\epsilon]$ denote
the ring of dual numbers. We call $T^1(Y,G):=\Def(Y,G)(k[\epsilon])$
the {\em tangent space} of the deformation functor $\Def(Y,G)$.
Theorem \ref{defthm} says in particular that there is a canonical
isomorphism
\begin{equation} 
  T^1(Y,G) \;\cong\; \HExt_G^1(\Ll_{Y/k},\Oo_Y).
\end{equation}
Moreover, Theorem \ref{defthm} together with standard arguments (see
e.g.\ \cite{Schlessinger68} or \cite{Vistoli99}) implies:

\begin{thm}  \label{defdat2thm}
  Suppose that $n:=\dim_k\HExt_G^1(\Ll_{Y/k},\Oo_Y)$ is finite (this
  holds, for instance, if $X$ is projective). Then $Y$ admits a
  versal deformation over a ring of the form
  \[
       R\univ \;=\; W(k)[[t_1,\ldots,t_n]]/\gen{f_1,\ldots,f_m}.
  \]
  If, moreover, $\HExt_G^2(\Ll_{Y/k},\Oo_Y)=0$ then $\Def(Y,G)$ is
  formally smooth and we have $R\univ=W(k)[[t_1,\ldots,t_n]]$.
\end{thm}

Let $Y_R$ be an equivariant deformation of $Y$ over $R$. Then the
quotient schemes $Z_R:=Y_R/G$ and $X_R:=Y_R/G$ are deformations of
$Z$ and $X$, respectively. Let $\Def(X;\tau_j)$ denote the functor
which classifies deformations of the {\em marked} curve $(X;\tau_j\mid
j\in B\branch)$, i.e.\ deformations $X_R$ of $X$ together with
sections $\tau_{j,R}:\Spec R\to X_R$ lifting the points $\tau_j$.

We claim that the association $Y_R\mapsto
X_R:=Y_R/G$ gives rise to a morphism of deformation functors
\begin{equation} \label{functorhom}
      \Def(Y,G)  \;\To\; \Def(X;\tau_j\mid j\in B\branch).
\end{equation}
To prove the claim we have to endow the curve $X_R$ with sections
$\tau_{j,R}:\Spec R\to X_R$ lifting the branch points $\tau_j$, for
all $j\in B\branch$. This is obvious for $j\in B\tame$: the $G$-action
on $Y_R$ induces an action of $H$ on $Z_R$ such that $X_R=Z_R/H$ and
such that the map $Z_R\to X_R$ is a tame $H$-cover lifting $Z\to X$.
It follows that the branch locus of $Z_R\to X_R$ is the disjoint union
of sections $\tau_{j,R}:\Spec R\to X_R$ lifting the tame branch points
$\tau_j$ (for $j\in B\tame$). Now let $j\in B\wild$ and let $\xi\in Z$
be a point above the wild branch point $\tau_j$. Let
$\phi_1,\ldots,\phi_s$ be a basis of $V$. We can choose this basis in
such a way that $\phi_1$ has a simple pole in $\xi$ and that
$\phi_2,\ldots,\phi_s$ generate the kernel of the residue map ${\rm
  res}_\xi:V\to\FF_p$. If we further replace $\phi_i$ by a multiple of
itself then we may assume that $\phi_i=\diff u_i/u_i$, with $\ord_\xi
u_1=1$ and $\ord_\xi u_i=0$ for $i>1$. In a neighborhood of $\xi$, the
cover $Y\to Z$ is (locally at $\xi$) given by $s$ Kummer equations
$y_i^p=u_i$, see Construction \ref{Yconstr}. Hence the deformation
$Y_R\to Z_R$ of $Y\to Z$ is (locally at $\xi$) given by $s$ Kummer
equations $y_i^p=u_{i,R}$, where $u_{i,R}$ lifts $u_i$. The equation
$u_{1,R}=0$ defines a section $\xi_R:\Spec R\to Z_R$ which lifts the
point $\xi$. We define $\tau_{j,R}:\Spec R\to X_R$ to be the image of
$\xi_R$. Using the $H$-action, one also checks that the definition of
$\tau_{j,R}$ for $j\in B\wild$ agrees with the definition of
$\tau_{j,R}$ for $j\in B\tame$, in case that $j\in B\tame\cap B\wild$.
This proves the claim.

It is well known (see e.g.\ \cite{DelMum69}) that the tangent
space of the deformation functor $\Def(X;\tau_j)$ is given by
\[
   T^1(X;\tau_j\mid j\in B\branch) \;\cong\; 
          H^1(X,\T_X(-\sum_{j\in B\branch}\,\tau_j)).
\]
Here $\T_X$ is the sheaf of tangent vectors of $X$. Hence the morphism
\eqref{functorhom} induces a $k$-linear map
\begin{equation} \label{ExtH1hom}
   \HExt_G^1(\Ll_{Y/k},\Oo_Y) \;\To\; 
       H^1(X,\T_X(-\sum_{j\in B\branch}\,\tau_j)).
\end{equation}
In the next subsection we will analyze this map in more detail.

\subsection{Analysis of $\HExt_G^1(\Ll_{Y/S},\Oo_Y)$} \label{defdat3}

By Lemma \ref{defdat1lem} and the spectral sequence \eqref{GOY4eq4} we
have 
\begin{equation} \label{defdat3eq2} 
  \EXt_G^n(\Ll_{Y/k},\Oo_Y) \;=\; \EXt_{G_0}^n(\Ll_{Y/k},\Oo_Y)^H 
           \;=\; \EXt_Y^n(\Ll_{Y/k},\Oo_Y)^G.
\end{equation}
The analogous statement for $\HExt_G^n$ holds as well. In other words,
we do not have to worry about group cohomology. In the following, we
shall use this sort of argument over and over again, sometimes without
mentioning it explicitly.

A {\em $V$-derivation} is an $\FF_p$-linear function
$\theta:V\to k(Z)$. We say that $\theta$ is integral at a point
$\xi\in Z$ if $\theta(\phi)\in\Oo_{Z,\xi}$ for all $\phi\in V$. Let
$\M$ be the sheaf of integral $V$-derivations; for an open
subset $U\in Z$ the group $\Gamma(U,\M)$ of sections over $U$ is the
set of $V$-derivations $\theta$ which are integral at each point
$\xi\in U$. Obviously, $\M$ is a locally free $\Oo_Z$-module whose rank
is equal to $s=\dim_{\FF_p}V$. There is a natural $H$-action on $\M$,
i.e.\ a structure of $H$-$\Oo_Z$-module, such that $\M^H$ is the sheaf
of $H$-equivariant $V$-derivations on $X$. 

Let $\T_Z=\Hom_Z(\Omega_{Z/k},\Oo_Z)$ denote the sheaf of tangent
vectors on $Z$. We write $\T_Z(D):=\T_Z\otimes\Oo_Z(D)$ etc. There is a
natural injection of $H$-$\Oo_Z$-modules
\[
       \T_Z(D) \;\inj\; \M
\]
which sends a vector field $\theta$ to its restriction to $V$. From the
definition of the divisor $D$ it is clear that this is well defined
and that the quotient $\M/\T_Z(D)$ is torsion free.

\begin{lem} \label{defdat3lem}
  There is a natural isomorphism of $H$-$\Oo_Z$-modules
  \begin{equation} \label{defdat3eq3}
      \M \;\liso\; \Hom_{G_0}(\Omega_{Y/k},\Oo_Y) :=\T_Y^{G_0}.
  \end{equation}
  Furthermore, we have a short exact sequence of $H$-$\Oo_Z$-modules
  \begin{equation} \label{defdat3eq4}
      0 \;\To\; \T_Z(-D_\infty) \;\To\; \M \;\To\;
                     \EXt_{G_0}^1(\Ll_{Y/k},\Oo_Y) \;\To\; 0.
  \end{equation}
\end{lem}

\proof It suffices to prove this locally on $Z$. Hence we may
assume that $Z=\Spec A$ is affine and that there is a basis
$\phi_1,\ldots,\phi_s$ of $V$ such that $\phi_i=\diff u_i/u_i$ with
$u_i\in A$ and such that $u_i$ has at most simple zeros. By
construction we have 
\[
     Y=\Spec B,\quad B=C/I,
\]
where $C=A[y_1,\ldots,y_s]$ is the polynomial algebra over $A$ in $s$
variables (with $V$-grading such that $y_i\in C_{\phi_i}$) and $I$ is
generated by the polynomials $u_i-y_i^p$. One checks that the
$B$-module $I/I^2$ is free, with $G_0$-invariant generators
$[u_i-y_i^p]$.

The cotangent complex $\Ll_{Y/k}$ may be identified
with the complex of $G$-$\Oo_Y$-modules associated to the complex
$L:=(I/I^2\to\Omega_{C/k}\otimes B)$ of $V$-graded $B$-modules with
$H$-action. The differential of this complex sends the generator
$[u_i-y_i^p]$ to the $1$-form $\diff u_i$. It follows that
$\Omega_{B/k}=H^0(L)$ is the direct sum of the free $B$-module
generated by $\diff y_i$ and the torsion module 
\[
     (\Omega_{B/k})_{\rm tors} \;=\; 
         \frac{\Omega_{A/k}}{\gen{\diff u_i}} \,\otimes_A B.
\]
Let $\theta:V\to A$ be an everywhere integral $V$-derivation. It gives
rise to a $G_0$-equivariant derivation $\eta:\Omega_{B/k}\to B$ which
is zero on $(\Omega_{B/k})_{\rm tors}$ and such that
\[
      \eta(\diff y_i) \;:=\; \theta(\phi_i)\,y_i.
\]
One checks that the association $\theta\mapsto \eta$ defines an
isomorphism of $H$-$\Oo_Z$-modules \eqref{defdat3eq3}.

Since both nontrivial terms of the complex $L$ are locally free
$B$-modules, we have
\[
    \HExt_{G_0}^n(\Ll_{Y/k},\Oo_Y) \;=\; H^n(\Hom_{G_0}^\bullet(L,B)).
\]
For $n=1$ this gives the exact sequence
\begin{equation}\label{defdat3eq5}
   \Hom_{G_0}(\Omega_{C/k},B) \;\To\; \Hom_{G_0}(I/I^2,B) \;\To\;
        \HExt_{G_0}^1(\Ll_{Y/k},\Oo_Y) \;\to\; 0.
\end{equation}
Let $\theta:V\to A$ be a global section of $\M$. There exists a unique
$B$-linear and $G_0$-equivariant map $\nu:I/I^2\to B$ such that
\[
        \nu([u_i-y_i^p]) \;=\; u_i\,\theta(\phi_i)
\]
for all $i$. This defines an $A$-linear map
\begin{equation} \label{defdat3eq6}
    H^0(Z,\M) \;\To\; \Hom_{G_0}(I/I^2,B).
\end{equation}
From \eqref{defdat3eq5} and \eqref{defdat3eq6} we obtain the sequence
\eqref{defdat3eq4}. It is easy to see that this sequence is
$H$-equivariant and does not depend on the choice of the basis of $V$.
It remains to show that \eqref{defdat3eq4} is exact.

Exactness on the left is obvious; exactness in the middle follows
easily from the exactness of \eqref{defdat3eq5}. To prove exactness on
the right, let $\nu:I/I^2\to B$ be a $B$-linear and $G_0$-equivariant
homomorphism.  We can define a $V$-derivation $\theta:V\to k(Z)={\rm
  Frac}(A)$ by setting
\[
      \theta(\phi_i) \:=\; \frac{\nu([u_i-y_i^p])}{u_i}.
\]
By construction the images of $\theta$ and of $\nu$ in
$\HExt_{G_0}^1(\Ll_{Y/k},\Oo_Y)$ are equal. If $\theta$ was integral
everywhere then we would be done. However, if $\xi\in Z$ is a point in
the support of $D_\infty$ then $\theta$ may not be integral at $\xi$.
In this case we may suppose that $\ord_\xi u_1=1$ and that $\ord_\xi
u_i=0$ for $i>1$.  After shrinking $Z$ to a sufficiently small
neighborhood of $\xi$ we may suppose that
\[
       \theta' \;=\; \theta - \nu([u_1-y_1^p])\,\partial/\partial u_1|_V
\]
is integral everywhere. But since \eqref{defdat3eq5} is exact in the
middle, the image of $\theta'$ in $\HExt_{G_0}^1(\Ll_{Y/k},\Oo_Y)$ is
the same as the image of $\theta$. This finishes the proof of the
lemma.  \Endproof

Applying the exact functor $\F\mapsto\F^H$ to \eqref{defdat3eq4} we
obtain a short exact sequence of $\Oo_X$-modules
\begin{equation} \label{defdat3eq7}
    0 \;\To\; \T_Z(-D_\infty)^H \;\To\; \M^H \;\To\;
                     \EXt_G^1(\Ll_{Y/k},\Oo_Y) \;\To\; 0.
\end{equation}
The following proposition identifies the first boundary map associated
to \eqref{defdat3eq7} with the differential of the morphism of
deformation functors \eqref{functorhom}.

\begin{prop} \label{defdat3prop}
  The following diagram commutes:
  \begin{equation} \label{defdat3diag1}
  \begin{CD}
     \HExt_G^1(\Ll_{Y/k},\Oo_Y) @>>> H^0(X,\EXt_G^1(\Ll_{Y/k},\Oo_Y)) \\
     @VVV                         @VV{\partial}V                  \\
     H^1(X,\T_X(-\sum_{j\in B\branch}\,\tau_j)) @>{\cong}>>
       H^1(X,\T_Z(-D_\infty)^H) \\
  \end{CD}
  \end{equation}
  Here the upper horizontal arrow is deduced from the local-global
  spectral sequence \eqref{GOY5eq3}. The left vertical arrow is the
  tangent map of the morphism \eqref{functorhom}. The right vertical
  arrow is the boundary map of the short exact sequence
  \eqref{defdat3eq7}. The lower horizontal arrow comes from the
  canonical isomorphism
  \[
       \T_Z(-D_\infty)^H \;\cong\; \T_X(-\sum_{j\in
            B\branch}\,\tau_j).
  \]
\end{prop}

\proof 
Let us denote by $\Def(Z,H,D_\infty)$ the functor which classifies
$H$-equivariant deformations of the marked curve $(Z,D_\infty)$
(here we identify $D_\infty$ with its support, which consists of the
points of $Z$ lying above the wild branch points). By
\cite{BertinMezard00} the tangent space of $\Def(Z,H,D_\infty)$ is
canonically isomorphic to 
\[
    H^1(Z,\T_Z(-D_\infty)^H) \;\cong\; 
        H^1(X,\T_X(\,-\!\!\sum_{j\in B\branch }\tau_j)).
\]
Using this fact, the proposition is easily reduced to the case $H=1$
and $Z=X$.  

Let $Y'$ be a $G$-equivariant deformation of $Y$ over $R=k[\epsilon]$
and $Z':=Y'/G$ the induced deformation of $Z$. We have seen in the
last subsection that $Z'$ is naturally endowed with a lift $D_\infty'$
of the divisor $D_\infty$. We denote by $e(Y')$ the global section of
$\EXt_{G_0}^1(\Ll_{Y/k},\Oo_Y)$ corresponding to $Y'$, see \S \ref{ext}
and \S \ref{def}. Similarly, we denote by $e(Z',D_\infty')\in
H^1(Z,\T_Z(-D_\infty))$ the cohomology class representing
$(Z',D_\infty')$. We have to show that $e(Z')$ is the image of $e(Y')$
under the boundary map $\partial$.

To prove this, we will first recall the definition of $e(Y')$ and
$e(Z',D_\infty')$. Let $(U_\mu)$ be a covering of $Z$ by sufficiently
small affine open subsets $U_\mu=\Spec A_\mu$. Let $W_\mu=\Spec
B_\mu\subset Y$ be the inverse image of $U_\mu$. Also, let
$U_\mu'=\Spec A_\mu'\subset Z'$ (resp.\ $W_\mu'=\Spec B_\mu'\subset
Y'$) be the induced deformation of $U$ (resp.\ the induced
$G$-equivariant deformation of $W$).

Since $Z$ is smooth over $k$ there exists, for all $\mu$, a
(non-canonical) isomorphism of $R$-algebras
\[
    \sigma_\mu:\,A_\mu' \;\liso\; A_\mu\otimes_k R
\]
which lifts the identity on $A_\mu$. For each pair of indices
$\mu,\lambda$ we set $U_{\mu,\lambda}:=U_\mu\cap U_\lambda=\Spec
A_{\mu,\lambda}$. Then the equality
\[
    \sigma_\lambda\circ\sigma_\mu^{-1} \;=\;
       \Id_{A_{\mu,\lambda}} \;+\; \epsilon\cdot\theta_{\mu,\lambda}
\]
defines a vector field
$\theta_{\mu,\lambda}\in\Gamma(U_{\mu,\lambda},\T_Z)$. The $1$-cocycle
$(\theta_{\mu,\lambda})$ represents the cohomology class
$e(Z',D_\infty')$.

We may assume that
\[
   B_\mu \;=\; A[\,y_i \mid u_{\mu,i}-y_{\mu,i}^p \,],
\]
with $u_{\mu,i}\in A_\mu$ such that $\phi_i=\diff
u_{\mu,i}/u_{\mu,i}$. There is a $V$-grading on $B_\mu$ such that
$y_{\mu,i}\in(B_\mu)_{\phi_i}$. Let $y_{\mu,i}'\in(B_\mu')_{\phi_i}$
be a lift of $y_{\mu,i}$. Then $u_{\mu,i}':=(y_{\mu,i}')^p\in
A_\mu$. Note that $u_{\mu,i}'$ is independent of the choice of the
lift $y_{\mu,i}'$. Set
\[
     \sigma_\mu(u_{\mu,i}') \;=\; u_{\mu,i} + \epsilon\cdot v_{\mu,i}.
\]
Let $\nu_\mu$ be the section of the sheaf $\M$ over $U_\mu$ such that
\[
          \nu_\mu(\phi_i) \;=\; v_{\mu,i}.
\]
Using the definition of $e(Y')$ via Theorem \ref{defthm}, together
with the proof of Proposition \ref{defdat3prop}, one checks that the
image of $\nu_\mu$ under the second map in \eqref{defdat3eq6} is equal
to the restriction of $e(Y')$ to $U_\mu$. 

A straightforward computation shows that
\[
    \theta_{\mu,\lambda}(\phi_i) \;=\; 
      \frac{\theta_{\mu,\lambda}(u_{\lambda,i})}{u_{\lambda,i}}
        \;=\; v_{\mu,i} - v_{\lambda,i}
\]
for all $\mu,\lambda,i$. This means that $\theta_{\mu,\lambda}$ is
mapped to $\nu_\mu-\nu_\lambda\in\Gamma(U_{\mu,\lambda},\M)$ under the
first map in \eqref{defdat3eq6}. Therefore, $e(Z',D_\infty')$ is the
image of $e(Y')$ under the boundary map $\partial$. This is what we
wanted to prove.  \Endproof

\begin{thm} \label{defdat3thm}
  Suppose $H^n(X,\M^H)=0$ for $n=0,1$. Then the morphism
  \[
     \Def(Y,G) \;\To\; \Def(X;\tau_j)
  \]
  is an isomorphism. In particular, the deformation functor
  $\Def(Y,G)$ is unobstructed. 
\end{thm}

\proof
The hypothesis implies that the boundary map 
\begin{equation} \label{defdat3eq8}
    \partial:H^0(X,\EXt_G^1(\Ll_{Y/k},\Oo_Y)) \;\liso\; 
                      H^1(X,\T_Z(-D_\infty)^H)
\end{equation}
deduced from the short exact sequence \eqref{defdat3eq7} is an
isomorphism. The local-global spectral sequence for $\HExt_G^n$ gives
rise to a short exact sequence
\[
   0 \;\to\; H^1(X,\T_Y^G) \;\To\; \HExt_G^1(\Ll_{Y/k},\Oo_Y)
     \;\To\; H^0(X,\EXt_G^1(\Ll_{Y/k},\Oo_Y)) \;\to\; 0.
\]
But \eqref{defdat3eq3} and the hypothesis show that $H^1(X,\T_Y^G)=0$.
Therefore, it follows from Proposition \ref{defdat3prop} that the morphism
$\Def(Y,G)\to\Def(X;\tau_j)$ induces an isomorphism on tangent spaces. The
theorem would follow if we knew that $\Def(Y,G)$ is unobstructed.

The local global spectral sequence for $\HExt_G^n$ also shows that
\begin{equation} \label{defdat3eq9}
    \HExt_G^2(\Ll_{Y/k},\Oo_Y) \;=\; H^1(X,\EXt_G^1(\Ll_{Y/k},\Oo_Y)).
\end{equation}
Using again the long exact cohomology sequence deduced from
\eqref{defdat3eq8} and the hypothesis we see that \eqref{defdat3eq7}
is zero. Hence $\Def(Y,G)$ is unobstructed by Theorem
\ref{defthm}. This concludes the proof of the theorem. 
\Endproof

\begin{rem} 
  Suppose that all elements of $V$ are regular, i.e.\ 
  $B\wild=\emptyset$. Then we may regard $V$ as an $\FF_p$-subvector
  space of the $\chi$-isotypical part of $J_Z[p](k)$. It can be shown
  that the hypothesis $H^n(X,\M^H)=0$ of Theorem \ref{defdat3thm} is
  equivalent to the condition that the $\chi$-isotypical part of the
  group scheme $J_Z[p]$ is \'etale. Using this fact one can give a
  different proof of Theorem \ref{defdat3thm}. In the special case
  $\dim_{\FF_p}V=1$ this is the approach taken in \cite{special}.
\end{rem}


\section{Special deformation data} \label{special}

In this section we suppose that $X=\PP^1_k$. We begin by defining a
certain class of multiplicative deformation data over $X$, which we
call {\em special}. The definition of speciality may seem a little bit
ad hoc. However, we show that the deformation functor $\Def(Y,G)$
associated to $(Z,V)$ in the last section has some very nice
properties if $(Z,V)$ is special. These are the {\em lifting property}
(Theorem \ref{special2thm}), the {\em local-global principle} (Theorem
\ref{special3thm}) and {\em rigidity} (Theorem \ref{special4thm}). At
a deeper level, these properties are explained by the way special
deformation data arise in the study of three point covers with bad
reduction, see \cite{special} and \cite{bad}.

We also prove a technical result (Proposition \ref{special5prop})
which is used in \cite{bad}.

\subsection{}  \label{special1}

Let $p$ be a prime, $H$ a finite group of order prime to $p$ and
$\chi_0:H\to\bar{\FF}_p^\times$ a one dimensional character on $H$
with values in the algebraic closure of $\FF_p$. The values of
$\chi_0$ generate a finite field $\FF_q$ with $q=p^s$ elements. Set
\[
    \chi_i \;:=\; \chi_0^{p^i}, \qquad
      \chi \;=\; \sum_{i=0}^{s-1}\, \chi_i.
\]
Then $\chi$ is an irreducible $\FF_p$-valued character. 

Let $k$ be an algebraically closed field of characteristic $p$ and set
$X:=\PP^1_k$.  Let $(Z,V)$ be a (multiplicative) deformation datum of
type $(H,\chi)$ over $X$. Choose a basis
$\omega_0,\ldots,\omega_{s-1}$ of $V\otimes_{\FF_p}\bar{\FF}_p$
consisting of eigenvectors, such that
\begin{equation} \label{special1eq1}
  \alpha^*\omega_i \;=\; \chi_i(\alpha)\,\omega_i,
\end{equation}
for all $\alpha\in H$ and $i=0,\ldots,s-1$. Let $\mathcal{C}$ denote
the {\em Cartier operator}. Since $\mathcal{C}$ is $p^{-1}$-linear and
is the identity on $V$, we have
\begin{equation} \label{special1eq2}
  \mathcal{C}(\omega_{i+1}) \;=\; c_i\,\omega_i
\end{equation}
for a constant $c_i\not=0$ in $\bar{\FF}_p$.  (Here and for the rest
of this section we will consider the index $i$ modulo $s$.) After
multiplying the $\omega_i$ with a constant in $\bar{\FF}_p$ we may
assume that $c_i=1$. 
 
As in the previous section, we denote by $\tau_j$, $j\in B$, the
critical points on $X$. Choose $j\in B$ and a point $\xi\in Z$
above $\tau_j$, and set
\[
    m_j \;:=\; |{\rm Stab}_H(\xi)|, \qquad
    h_j^{(i)} \;:=\; \ord_\xi\omega_i + 1, \qquad
    \sigma_j^{(i)} \;:=\; \frac{h_j^{(i)}}{m_j}.
\]
The tuple $(\sigma_j^{(i)})_{i,j}$ is called the {\em signature} of
the deformation datum $(Z,V)$. For a rational number $w$, we let
$\gen{w}$ denote the fractional part of $w$ (such that
$0\leq\gen{w}<1$ and $w-\gen{w}\in\ZZ$). 

\begin{lem} \label{special1lem}
For all $i$ we have
\begin{enumerate}
\item
  $\qquad\sum_{j\in B}\,(\sigma_j^{(i)}-1) \;=\; -2$,
\item
  $\qquad\gen{\sigma_j^{(i)}} \;=\; \gen{p^i\sigma_j^{(0)}}$.
\end{enumerate}
\end{lem}   
  
\proof Part (i) follows from a straightforward computation using the
Riemann--Hurwitz formula. To prove (ii), let $\xi\in Z$ be a point
above $\tau_j$ and $z$ a local coordinate at $\xi$. The inertia
character $\psi_\xi:H_\xi\to k^\times$ is determined by the congruence
\begin{equation} \label{special1eq4}
  \alpha^*z \;\equiv\; \psi_\xi(\alpha)\,z \pmod{z^2}
\end{equation}
for all $\alpha\in H_\xi:={\rm Stab}_H(\xi)$.  Now \eqref{special1eq1}
and the definition of $h_j^{(i)}$ imply that
\[
     \chi_0^{p^i} \;=\; \psi_\xi^{h_j^{(i)}}.
\]
Part (ii) of the lemma follows.
\Endproof 

\begin{defn}
  The deformation datum $(Z,V)$ is called {\em pure} if for all $i$
  we have
  \[
     \sum_{j\in B}\; \gen{\sigma_j^{(i)}} \;=\; 1.
  \]
\end{defn}

\begin{lem} \label{special1lem2}
  Let $\M$ be the sheaf of $H$-$\Oo_Z$-modules defined in \S
  \ref{defdat3}.  The deformation datum $(Z,V)$ is
  pure if and only if $H^n(X,\M^H)=0$ for $n=0,1$. 
\end{lem} 
       
\proof 
Let $\Oo_{Z,\chi_i}$ denote the
$\chi_i$-isotypical part of the sheaf $\pi_*\Oo_Z$. 
Let $\theta:V\to k(Z)$ be an $H$-equivariant $V$-derivation and extend
it $k$-linearly to $V_k$. Then $f_i:=\theta(\omega_i)$ is a
meromorphic section of $\Oo_{Z,\chi_i}$; it is holomorphic at $\tau\in
X$ if and only if $\theta$ is integral at all points $\xi\in Z$ above
$\tau$. Therefore, the rule $\theta\mapsto (f_i)$ defines an
isomorphism
\[
      \M^H \;\liso\; \bigoplus_{i=0}^{s-1}\; \Oo_{Z,\chi_i}
\]
of $\Oo_X$-modules. A local calculation
as in the proof of Lemma \ref{special1lem} shows that 
\[
    \deg\,\Oo_{Z,\chi_i} \;=\; 
        - \sum_{j\in B}\, \gen{\sigma_j^{(i)}}.
\]
Hence the lemma follows from the Riemann--Roch formula.
\Endproof

\subsection{}  \label{special2}

In order to discuss special deformation data we need some more notation:

\begin{notation} \label{special2not}
  Let $(Z,V)$ be a deformation datum of type $(H,\chi)$ and signature
  $(\sigma_j^{(i)})$. Set 
  \[
     \nu_j^{(i)} \;:=\; \lfloor\sigma_j^{(i)}\rfloor \;=\; 
         \sigma_j^{(i)} - \gen{\sigma_j^{(i)}}, \qquad
     \nu_j \;:=\; \min_i \,\nu_j^{(i)},
  \]
  and
  \[
     a_j^{(i)} \;:=\; m_j\cdot\gen{\sigma_j^{(i)}}, \qquad
       a_j \;:=\; \min_i\,a_j^{(i)}.
  \]
\end{notation}

\begin{defn} \label{special2defn}
  The deformation datum $(Z,V)$ is called {\em special} if 
  $\sigma_j^{(i)}\not=1$ for all $i$ and $j$ and if the
  following holds. There exists a subset $B_0\subset B$ with exactly
  three elements such that
  \[
      \nu_j \;=\; \begin{cases}  \quad 0,& \qquad j\in B_0\\
                                 \quad 1,& \qquad j\not\in B_0.\\
                  \end{cases}
  \]
  A special deformation datum $(Z,V)$ is called {\em normalized} if 
  $\{\,\tau_j\mid j\in B_0\,\}=\{0,1,\infty\}\subset X=\PP^1$.
\end{defn}

For the rest of this section we assume that the deformation datum
$(Z,V)$ is special. Whenever it is convenient, we may also assume that
$(Z,V)$ is normalized. (However, sometimes it more convenient to have
$\tau_j\not=\infty$ for all $j\in B$.) Since $\sigma_j^{(i)}\not= 1$ we
conclude that $B\wild\subset B_0$ and that $B\branch=B$. We set
\[
    B\new \;:=\; B-B_0, \qquad B\prim \;:=\; B_0-B\wild.
\]
For an explanation of the terminology, see \cite{bad}.

\begin{lem} \label{special2lem}
  Suppose $(Z,V)$ is special. Then the following holds.
  \begin{enumerate}
  \item
    The deformation datum $(Z,V)$ is pure. 
  \item
    We have $\nu_j^{(i)}=\nu_j$ for all $j\in B$.
  \item
         Let $j\in B-B\wild$ and let $\xi\in Z$ be a point above
         $\tau_j$. Then for all $\phi\in V$ we have
    \[
        \ord_\xi\phi \;=\; \nu_jm_j + a_j -1.
    \]
  \end{enumerate}
\end{lem}

\proof
By Lemma \ref{special1lem} (i) we have
\begin{equation} \label{special2eq3}
   1 \;\;=\;\;  3 \,+\, \sum_{j\in B}\,(\nu_j^{(i)}-1) 
              \,+\, \sum_{j\in B} \gen{\sigma_j^{(i)}} 
     \;\;\geq\;\;  \sum_{j\in B} \gen{\sigma_j^{(i)}}.
\end{equation}
Suppose that $\sum_j\gen{\sigma_j^{(i)}}=0$. Since
$\sigma_j^{(i)}\not=1$, it would follow that $B=B_0$ and
$\sigma_j^{(i)}=0$ for all $i$ and $j$. But then
$\sum_j(\sigma_j^{(i)}-1)=-3$, contradicting Lemma \ref{special1lem}.
We conclude that $\sum_j\gen{\sigma_j^{(i)}}=1$, proving (i). We also
conclude that the inequality \eqref{special2eq3} is an equality, which
means that $\nu_j^{(i)}=\nu_j$, proving (ii).   

It follows from Lemma \ref{special1lem} (ii) that
$\ord_\xi\omega_i=m_j\sigma_j^{(i)}-1$ takes pairwise distinct values
for all $i$. Every element $\phi\in V$ can be written as $\sum_i
c_i\omega_i$, with $c_i\in k$. Using $\mathcal{C}(\phi)=\phi$ and
Equations \eqref{special1eq1} and \eqref{special1eq2}, one shows that
$c_i\not=0$ for all $i$. Therefore,
\[
      \ord_\xi\phi \;=\; \min_i\, (\ord_\xi\omega_i) \;=\;
         \min_i\,(m_j\nu_j^{(i)} + a_j^{(i)} -1).
\]
Now (iii) follows from (ii). 
\Endproof

Putting Theorem \ref{defdat3thm}, Lemma \ref{special1lem2} and
Lemma \ref{special2lem} (i) together, we get:

\begin{thm}[Lifting property] \label{special2thm}
  Let $Y$ be the curve with $G$-action corresponding to the special
  deformation datum $(Z,V)$, as defined in \S \ref{defdat}. The
  homomorphism of deformation functors
  \[
       \Def(Y,G) \;\To\; \Def(X;\tau_j)
  \]
  is an isomorphism. 
\end{thm}

\begin{prob} \label{special2prob}
  Let $(H,\chi)$ be as in the beginning of this section. Let
  $(\sigma_j^{(i)})$ be a tuple of rational numbers (indexed by $j\in
  B$ and $i\in \ZZ/s$) such that the statements of Lemma
  \ref{special1lem} and of Definition \ref{special2defn} hold. 
  Furthermore, let $(\tau_j)_{j\in B}$ be a $B$-tuple of closed points
  of $X=\PP^1_k$. Does there exists a special deformation datum
  $(Z,V)$ of type $(H,\chi)$ with signature $(\sigma_j^{(i)})$ and
  critical points $(\tau_j)$?
\end{prob}

\begin{prop} \label{special2prop}
  With assumptions as in Problem \ref{special2prob}:
  \begin{enumerate}
  \item Suppose that the character $\chi_0:H\to\FF_q^\times$ is
    injective. Then if it exists, the special deformation datum
    $(Z,V)$ is uniquely determined (up to isomorphism) by the datum
    $(H,\chi,\sigma_j^{(i)},\tau_j)$.
  \item
    Fix $(H,\chi,\sigma_j^{(i)})$. The set of all tuples $(\tau_j)$
    such that there exists a special deformation datum $(Z,V)$ with
    critical points $(\tau_j)$ is a locally closed subset of
    $(\PP^1_{\FF_p})^B$. 
  \end{enumerate}
\end{prop}
We will see later (Theorem \ref{special4thm}) that the set of tuples
$(\tau_j)$ in (ii) is actually finite. 

\vspace{1ex} \proof (Compare with \cite{special}, \S 3.5) Suppose that
$(Z,V)$ exists. Let $\bar{H}:=H/\Ker(\chi_0)$,
$\bar{Z}:=Z/\Ker(\chi_0)$ and $\bar{\chi}$ the restriction of $\chi$
to $\bar{H}$. The subvector space $V\subset\Omega_{k(Z)/k}$ descends
to a subvector space $\bar{V}\subset\Omega_{k(\bar{Z})/k}$. One checks
that $(\bar{Z},\bar{V})$ is again a special deformation datum, of type
$(\bar{H},\bar{\chi})$. The signature $(\sigma_j^{(i)})$ and the set
$(\tau_j)$ of critical points remain unchanged during this descent.
Therefore, we may assume that $\chi_0$ is injective, even for the
proof of (ii). (For (i) the assumption of injectivity is necessary
because the cover $Z\to\bar{Z}$ is not unique if
$\Ker(\chi_0)\not=1$.)

If $\chi_0$ is injective then $H$ is cyclic of order $m$ where $m$ is
a positive integer such that $\FF_q=\FF_p[\zeta_m]$; in particular,
$m|q-1$. Set $a_j^{(i)}:=m\,\gen{\sigma_j^{(i)}}$. Then the $a_j^{(i)}$
are integers with $0\leq a_j^{(i)}<m$, $\sum_j\,a_j^{(i)}=m$ and
$a_j^{(i+1)} \equiv p^i\,a_j^{(i)} \pmod{m}$. The proof of Lemma
\ref{special1lem} shows that the $a_j^{(i)}$ determine the
ramification type of the $m$-cyclic cover $\pi:Z\to X=\PP^1_k$. This
can be made more explicit with Kummer theory. In fact, for all $i$
there exists a rational section $z_i$ of $\Oo_{Z,\chi_i}$ which
satisfies the equation
\begin{equation} \label{special2eq5}
       z_i^m \;=\; \prod_{j\in B} \,(x-\tau_j)^{a_j^{(i)}}.
\end{equation}
Here $x$ denotes the standard coordinate on $X=\PP^1$ and we assume,
without loss of generality, that $\tau_j\not=\infty$. The curve $Z$ is
the smooth projective model of the plane curve with equation
\eqref{special2eq5} (for any $i$). 

We claim that the eigenvector $\omega_i$ of $V_k$ is of the form
\begin{equation} \label{special2eq6}
   \omega_i \;=\; \epsilon_i\,
        \frac{z_i\,\diff x}{\prod_{j\in B_0}\,(x-\tau_j)}
\end{equation}
for some constant $\epsilon_i\in k$. Indeed, a local calculation shows
that the right hand side of \eqref{special2eq6} has everywhere the
right order of poles and zeros compatible with the signature
$(\sigma_j^{(i)})$ and the set of critical points $(\tau_j)$. This
proves the claim. If we plug in \eqref{special2eq6} into the equation
\begin{equation} \label{special2eq7}
      \mathcal{C}(\omega_{i+1}) \;=\; c_i\,\omega_i
\end{equation}
and look at Taylor series (say in $x$) on both sides, we obtain a set
of algebraic equations with coefficients in $\FF_p$ which are
satisfied by the tuple $(\tau_j)$. These equations define a Zariski
closed subset of $(\PP^1_{\FF_p})^B$. The conditions $c_i\not=0$
define an open subset of this closed subset. We have shown that the
set of tuples $(\tau_j)$ coming from a special deformation datum with
given type and signature are contained in a certain locally closed
subset of $(\PP^1_{\FF_p})^B$. It is clear that $(Z,V)$ is uniquely
determined by the datum $(H,\chi,\sigma_j^{(i)},\tau_j)$.

Conversely, let $(\tau_j)$ be a $B$-tuple of $k$-rational points of
$\PP^1$ which is contained in the locally closed subset constructed
above. This means that if we define an $H$-cover $\pi:Z\to X=\PP^1$ by
equation \eqref{special2eq5} and define differentials $\omega_i$ on
$Z$ by equation \eqref{special2eq6} then \eqref{special2eq7} holds
with certain constants $c_i\not=0$. After multiplying the $\omega_i$
by suitable constants we may assume that $c_i=1$. Let $V'\subset
\Omega_{k(Z)/k}$ be the $\FF_q$-linear subspace spanned by the
$\omega_i$. The Cartier operator $\mathcal{C}$ stabilizes $V'$ and
acts semi-simply on it. A well known lemma in $p^{-1}$-linear algebra
shows that the stabilizer $V$ of $\mathcal{C}$ inside $V'$ is an
$\FF_p$-vector space of dimension $s=\dim_{\FF_q}V'$. Here we can be
more explicit: if $\alpha\in H$ is an element such that
$\chi_0(\alpha)$ generates $\FF_q$ then
\begin{equation} \label{special2eq8}
      \phi_l \;:=\; \sum_i\; \chi_{i+l}(\alpha)\cdot\omega_i, \qquad
           l=0,\ldots,s-1
\end{equation}
gives a basis for $V$. By construction, $(Z,V)$ is a special
deformation datum of type $(H,\chi)$, signature $(\sigma_j^{(i)})$ and
with critical points $(\tau_j)$. This concludes the proof of the
proposition. \Endproof

\subsection{The local-global principle} \label{special3}

For $j\in B$, let $\Yd_j$ denote the completion of $Y$ at the critical
point $\tau_j$, see \S \ref{def3}. Given an equivariant deformation
$Y_R$ of $Y$, we denote by $\Yd_{j,R}$ the completion of $Y_R$ at
$\tau_j$; this is an equivariant deformation of $\Yd_j$. We obtain a
morphism
\[
     \Phi:\,\Def(Y,G) \;\To\; \prod_{j\in B}\,\Def(\Yd_j,G)
\]
which maps a deformation $Y_R$ to the tuple $(\Yd_{j,R})_j$. Following
\cite{BertinMezard00}, we call $\Phi$ the {\em local-global morphism}.
By the results of \S \ref{def3} we can identify the natural morphism
arising from the local-global spectral sequence
\begin{equation} \label{special3eq2}
    \HExt_G^1(\Ll_{Y/k},\Oo_Y) \;\To\; 
       \bigoplus_{j\in B}\; \EXt_G^1(\Ll_{Y/k},\Oo_Y)\sphat_{\tau_j} 
\end{equation}
with the differential of $\Phi$. In contrast to the situation studied
in \cite{BertinMezard00}, $\Phi$ is not formally smooth unless $s=1$.
In fact, if $s>1$ then the groups
$\EXt_G^1(\Ll_{Y/k},\Oo_Y)\sphat_{\tau_j}$ are not finite-dimensional
over $k$. However, if we restrict our attention to the image of
$\Phi$, then we obtain a {\em local-global-principle}, comparable to
\cite{BertinMezard00}, Th\'eor\`eme 3.3.4.
 
\begin{lem}  \label{special3lem}
  The map \eqref{special3eq2} is injective. Its image is the direct
  sum
  \[
     \bigoplus_{j\in B\new}\T_{X,\tau_j}\otimes k(\tau_j) \;\subset\;
       \bigoplus_{j\in B}\,\EXt_G^1(\Ll_{Y/k},\Oo_Y)\sphat_{\tau_j}. 
  \]
\end{lem}

\proof We have already seen in the proof of Theorem \ref{defdat3thm}
that the natural map
\[
   \HExt_G^1(\Ll_{Y/k},\Oo_Y) \;\To\; H^0(X,\EXt_G^1(\Ll_{Y/k},\Oo_Y))  
\]
is an isomorphism. (Note that the hypothesis of Theorem
\ref{defdat3thm} is verified by Lemma
\ref{special2lem} (i).) Furthermore, we have an isomorphism
\begin{equation} \label{special3eq3}
    H^0(X,\EXt_G^1(\Ll_{Y/k},\Oo_Y)) \;\cong\; 
      H^1(X,\T_X(\,-\!\sum_{j\in B}\tau_j)).
\end{equation}
The $k$-dimension of \eqref{special3eq3} is $|B|-3=|B\new|$ by 
Riemann-Roch. 

Let $\E_{\rm tor}\subset\EXt_G^1(\Ll_{Y/k},\Oo_Y)$ be the maximal
sub-$\Oo_X$-module which is torsion. The sequence
\eqref{defdat3eq7} and a local computation shows 
\[
   \E_{\rm tor} \;\cong\; \frac{\T_Z(D)^H}{\T_Z(-D_\infty)^H} 
      \;\cong\; \frac{\T_X(\;\sum_{j\in B}(1-\nu_j)\tau_j)}
                     {\T_X(\,-\sum_{j\in B}\tau_j)}.
\]
Therefore,
\begin{equation} \label{special3eq4}
   H^0(X,\E_{\rm tor}) \;\cong\; 
     \bigoplus_{j\in B\new}\; \T_{X,\tau_j}\otimes k(\tau_j).
\end{equation}
Comparing dimensions, we find that $H^0(X,\E_{\rm tor})\inj
H^0(X,\EXt_G^1(\Ll_{Y/k},\Oo_Y))$ is an isomorphism.  This proves the
lemma.  \Endproof

Let $\Def(\Yd_j,G)^\dagger\subset\Def(\Yd_j,G)$ denote the image of
$\Def(Y,G)$ under the localization map. In other words,
$\Def(\Yd_j,G)^\dagger$ classifies those equivariant deformations of
$\Yd_j$ which arise as the completion of a global deformation of $Y$.
We denote by 
\[
     \Phi^\dagger:\,\Def(Y,G) \;\To\; \Def(Y,G)^{\rm loc}:=
        \prod_{j\in B} \Def(\Yd_j,G)^\dagger
\]
the restriction of $\Phi$ onto its image. By Lemma \ref{special3lem},
the differential of $\Phi^\dagger$ is the isomorphism
\[
    \HExt_G^1(\Ll_{Y/S},\Oo_Y) \;\liso\; 
       \bigoplus_{j\in B\new}\; \T_{X,\tau_j}\otimes k(\tau_j).
\]

\begin{thm}[Local-global principle] \label{special3thm}\ 
\begin{enumerate}
\item
  The functor $\Def(\Yd_j,G)^\dagger$ admits a versal
  deformation  over the ring
  \[
     \Rt_j \;:=\quad \begin{cases}
               \quad W(k), &\qquad \text{\rm for $j\in B_0$,} \\
        \quad W(k)[[t_j]], &\qquad \text{\rm for $j\in B\new$.}
                  \end{cases}
  \]
\item
  The functor $\Def(Y,G)$ admits an effective universal deformation.
  Let $\Rt$ be the universal deformation ring.
\item The restricted local-global morphism $\Phi^\dagger$ is an
  isomorphism. Therefore, we have
  \[
       \Rt \;\cong\; \hat{\otimes}_{W(k)}\,\Rt_j \;=\;
           W(k)[[\,t_j\mid j\in B\new\,]].
  \]
\end{enumerate}
\end{thm}

\proof By Theorem \ref{defdat2thm} the functor $\Def(Y,G)$ admits a versal deformation. By Theorem \ref{defdat3thm} it is
unobstructed. The space of `infinitesimal automorphisms' of
$\Def(Y,G)$ is isomorphic to $H^0(X,\T_Y^G)$. By Lemma
\ref{defdat3lem} we have
\[
    H^0(X,\T_Y^G) \;\cong\; H^0(X,\M^H) \;=\; 0.
\]
Therefore, the versal deformation of $\Def(Y,G)$ is also universal,
see \cite{Schlessinger68}. It is effective by Grothendieck's existence
theorem. This proves (ii).  The functors
$\Def(\Yd_j,G)^\dagger$ are unobstructed, because the same holds for
$\Def(Y,G)$. Using Theorem \ref{def3thm} and Lemma \ref{special3lem},
it is easy to verify Schlessinger's axioms \cite{Schlessinger68},
showing that $\Def(\Yd_j,G)^\dagger$ admits a versal deformation
over the ring $\Rt_j$, as claimed in (i).  Finally, Lemma
\ref{special3lem} (i) together with the argument used in the proof of
\cite{BertinMezard00}, Th\'eor\`eme 3.3.4, shows that $\Phi^\dagger$
is an isomorphism.  This finishes the proof of the theorem. \Endproof

\begin{rem}
  If $s=\dim_{\FF_p}V=1$ then $\Phi^\dagger=\Phi$. 
\end{rem}

\subsection{Rigidity} \label{special4}

In this subsection $(Z,V)$ and $Y$ will be as before. Let $R$ be an
Artinian local $k$-algebra with residue field $k$. Since $R$ has
characterisitc $p$, an equivariant deformation of $Y$ over $R$
corresponds to a deformation datum $(Z_R,V_R)$ over $R$ which lifts
$(Z,V)$. By this we mean that $\pi_R:Z_R\to X_R=\PP^1_R$ is a tamely
ramified $H$-Galois cover lifting $\pi:Z\to X$ and $V_R\subset
H^0(Z_R,\Omega_{Z_R/R}(D_{\infty,R}))$ is an $H$-stable $\FF_p$-vector
space of logarithmic differentials lifting $V$ (here
$D_{\infty,R}\subset Z_R$ is a relative Cartier divisor lifting
$D_\infty$).

Let $Y_R$ be an equivariant deformation of $Y$ and $(Z_R,V_R)$ the
corresponding deformation datum. Choose $j\in B\new$ and a point $\xi\in Z$
lying above $\tau_j$. By the theory of tame ramification, there exists a local
parameter $z$ for $Z_R$ at $\xi$ such that $\Oo_{Z_R,\xi}=R[[z]]$ and
$\alpha^*z=\psi(\alpha)\cdot z$ for some character $\psi:H_\xi\to R^\times$.
We say that the deformation $Y_R$ is {\em $j$-special} if every element
$\phi\in V_R$ is of the form
\[
    \phi \;=\; z^{m_j+a_j-1}(c_0+c_1z+\ldots)\,\diff z
\]
with $c_i\in R$ and $c_0\in R^\times$. Note that this condition is independent
of the choice of $z$. 

\begin{lem} \label{special4lem}
  The equivariant deformation $Y_R$ is trivial (i.e.\ isomorphic to
  $Y\otimes_k R$) if and only if it is $j$-special for all $j\in
  B\new$.
\end{lem}

\proof One direction of the claim follows immediately from Lemma
\ref{special2lem} (iii). To prove the other direction, suppose that
$Y_R$ is $j$-special for every $j\in B\new$. We have to show that
$Y_R$ is the trivial deformation. By Theorem \ref{special3thm} (iii)
it suffices to show that the completion $\Yd_{j,R}$ of $Y_R$ at
$\tau_j$ is the trivial deformation of $\Yd_j$, for all $j\in B\new$.
Fix one index $j$. Since $R$ is Artinian, we may prove triviality of
$\Yd_{j,R}$ by induction: suppose that $\Yd_{j,R'}:=\Yd_{j,R}\otimes_R
R'$ is trivial, where $R':=R/\m_R^n$ for some $n\geq 1$. Then we want
to conclude that $\Yd_{j,R''}$ is trivial, where $R'':=R/\m_R^{n+1}$.
To simplify the notation, we may even assume that $R=R''$.

The `difference' between $\Yd_{j,R}$ and the trivial deformation
$\Yd_j\otimes_k R$, considered as lifts of the trivial deformation
$\Yd_j\otimes_k R'$, is measured by an element $\bar{\theta}_j$ in
\[
    \EXt_G^1(\Ll_{Y/S},\Oo_Y)_{\tau_j}\sphat\otimes_k\m_R^n \;\cong\;
        (\M^H/\T_Z(-D_\infty)^H)_{\tau_j}\sphat\otimes_k\m_R^n,
\]
see Theorem \ref{def3thm} (ii). Since $\Yd_{j,R}$ lies in the image of the
local-global morphism $\Phi$ the element $\bar{\theta}_j$ lies in the subspace
\[
   (\T_Z(D)^H/\T_Z(-D_\infty)^H)_{\tau_j}\sphat\otimes_k\m_R^n \;\cong\;
     \T_{X,\tau_j}\otimes k(\tau_j)\otimes_k\m_R^n,
\]
see the proof of Lemma \ref{special3lem}. In other words, we may regard
$\bar{\theta}_j$ as a tangent vector at $\tau_j$, with values in the
$k$-vector space $\m_R^n$. We have to show that $\bar{\theta}_j=0$.

Choose a point $\xi\in Z$ above $\tau_j$. Let $z$ be a local parameter
of $Z_R$ at $\xi$ such that $\Od_{Z_R,\xi}=R[[z]]$ and
$\alpha^*z=\psi_\xi(\alpha)\cdot z$ for a character $\psi_\xi:H_\xi\to
R^\times$. It follows that $\Od_{X_R,\tau_j}=R[[x]]$, where
$x:=z^{m_j}$ and $m_j:=|H_\xi|$. Note that the fiber product
$\Yd_{\xi,R}:=Y_R\times_{Z_R}\Spec R[[z]]$ is a connected component of
$\Yd_{j,R}$.  Let $\phi_1,\ldots,\phi_s$ be a basis of $V_R$. We have
$\phi_i=\diff u_i/u_i$ for a unit $u_i\in R[[z]]^\times$ which is
unique up to multiplication by a $p$th power.  The $G_0$-torsor
$\Yd_{\xi,R}\to\Spec R[[z]]$ is given by the Kummer equations
\[
     y_i^p \;=\; u_i, \qquad i=1,\ldots,s.
\]
By our induction hypothesis, the induced deformation $\Yd_{j,R'}$ is
trivial. This means that, for a suitable choice of the parameter $z$
and the units $u_i$, the image of $u_i$ in $R'[[z]]$ actually lies in
the subalgebra $k[[z]]\subset R'[[z]]$. In other words, we have
\[
    u_i \;=\; \ub_i + v_i, \qquad \ub\in k[[z]]^\times,\qquad
        v_i\in k[[z]]\otimes_k\m_R^n.
\]

{\bf Claim:} The tangent vector $\bar{\theta}_j$ extends to a vector
field $\theta_j\in\T_{X,\tau_j}\otimes\m_R^n$ in a neighborhood of
$\tau_j$, with values in $\m_R^n$, such that
\begin{equation} \label{special4eq2}
       v_i \;=\; \theta_j(\diff\ub_i), 
\end{equation}
for $i=1,\ldots,s$.

Let us prove this claim. The class in
$\EXt_G^1(\Ll_{Y/S},\Oo_Y)_{\tau_j}\sphat\otimes_k\m_R^n$ corresponding
to $\bar{\theta}_j$ lifts to a local section $\theta_j'$ of the sheaf
$\M^H$ in a neighborhood of $\tau_j$, via the exact sequence
\eqref{defdat3eq7}. We consider $\theta_j'$ as an $\FF_p$-linear and
$H$-equivariant map $\theta_j:V\to k[[z]]\otimes\m_R^n$. By the
definition of $\bar{\theta}_j$ in terms of the deformation
$\Yd_{j,R}$, we have $\theta_j'(\phib_i)=v_i/\ub_i$ (compare with the
proof of Proposition \ref{defdat3prop}). But since $\Yd_{j,R}$ lies in
the image of the local-global morphism, $\theta_j'$ is actually the
restriction to $V$ of a vector field $\theta_j$ on $X$ which is
regular in a neighborhood of $\tau_j$ (compare with the proof of Lemma
\ref{special3lem}). The claim follows.

We can now finish the proof of the lemma. The vector field $\theta_j$
appearing in the claim we have just proved can be written as follows:
\[
   \theta_j \;=\; (b_0+b_1x+\ldots)\,\vf{x} \;=\; 
     \frac{1}{m_j}(b_0z^{1-m_j}+b_1z+\ldots)\,\vf{z},
\]
with $b_\mu\in\m_R^n$. Since by assumption the deformation $Y_R$ is
$j$-special we have $\diff u_i=z^{m_j+a_j-1}(c_0+c_1z+\ldots)\,\diff
z$ with $\cb_0\not=0$. From \eqref{special4eq2} we get 
\begin{equation} \label{special4eq3}
     v_i \;=\; 
        \frac{z^{a_j}}{m_j}\;(\,\cb_0b_0\,+\,
                    (\cb_0b_1+\cb_1b_0)z^{m_j}\,+\,\ldots).
\end{equation}
But since $\diff u_i=\diff\ub_i+\diff v_i$ is divisible by
$z^{m_j+a_j-1}$ it follows that 
\[
      a_j\,\cb_0 \,b_0 \;=\;0.
\]
But $a_j$ is prime to $p$ and $\cb_0\not=0$, hence $b_0=0$.  We conclude
that $\bar{\theta}_j=0$, which completes the proof of the lemma.
\Endproof

\begin{thm}[Rigidity] \label{special4thm}
  There exist, up to isomorphism, at most a finite number of special
  deformation data of given type $(H,\chi)$. Moreover, every special
  deformation datum can be defined over a finite field.
\end{thm}

\proof For a fixed type $(H,\chi)$ there exists at most a finite
number of possibilities for the signature $(\sigma_j^{(i)})$ of a
special deformation datum. Therefore, we may also fix the signature
$(\sigma_j^{(i)})$. Let $U\subset(\PP^1)^B$ be the locally closed
subset from Proposition \ref{special2prop} (ii). Let
$U'\subset(\PP^1)^{B\new}$ denote the intersection of $U$ with the
closed subset $\{(0,1,\infty)\}\times(\PP^1)^{B\new}\subset(\PP^1)^B$.
Thus, a point on $U'$ corresponds to the branch locus of a {\em
  normalized} special deformation datum. To prove the corollary, it
suffices to show that $U'$ has pure dimension $0$. Suppose that $U'$
has an irreducible component of dimension $>0$. Then there exists an
algebraically closed field $k$ of characteristic $p$ and a nonconstant
morphism $\varphi:\Spec R\to U'$, with $R=k[[t]]$. We will show that
$\varphi$ is constant, which gives a contradiction.

Going again through the proof of Proposition \ref{special2prop}, we
see that $\varphi$ corresponds to a deformation datum $(Z_R,V_R)$
defined over $R$.  Moreover, the special and the generic fiber of
$(Z_R,V_R)$ are special deformation data. Applying Lemma
\ref{special2lem} to the generic fiber of $(Z_R,V_R)$, we see that the
pullback $(Z_{R'},V_{R'})$ of $(Z_R,V_R)$ over $R':=R/t^n$ is a
deformation of its special fiber which is $j$-special, for all $j\in
B\new$ and for all $n$. Hence it follows from Lemma \ref{special4lem}
that the curve $Y_{R'}$ corresponding to $(Z_{R'},V_{R'})$ is the
trivial deformation of its special fiber. By Theorem
\ref{special2thm}, this implies that the branch locus of the induced
$G$-cover $Y_{R'}\to X_{R'}=\PP^1_{R'}$ is constant, for all $n$. We
conclude that $\varphi:\Spec R\to U'$ is constant, which proves the
theorem.  \Endproof

\subsection{}  \label{special5}

In this last section we prove a proposition which links two of our
previous results on special deformation data: the lifting property
(Theorem \ref{special2thm}) and the local-global principle (Theorem
\ref{special3thm}). This proposition is a key ingredient for the proof
of the main result of \cite{bad}.

Let $\Y$ be the universal equivariant deformation of $Y$ over $\Rt$,
see Theorem \ref{special3thm}. The quotient scheme $\X:=\Y/G$ is
naturally equipped with sections $\tau_{j,\Rt}:\Spec\Rt\to\X$ lifting
the critical points $\tau_j$. We may suppose that $\X=\PP^1_{\Rt}$ and
that $\{\,\tau_{j,\Rt}\mid j\in B_0\,\}=\{0,1,\infty\}$. With this
normalization, we may regard the sections $\tau_{j,\Rt}$ for $j\in
B\new$ simply as elements of the ring $\Rt$. Let $[\tau_j]\in W(k)$
denote the Teichm\"uller lift of $\tau_j\in k$ and set
$T_j:=\tau_{j,\Rt}-[\tau_j]$. By the lifting property (Theorem
\ref{special2thm}) we have
\[
   \Rt \;=\; W(k)[[\,t_j\mid j\in B\new\,]] \;=\; 
                W(k)[[\,T_j\mid j\in B\new\,]].
\]
A priori, it is not clear that these two sets of coordinates of $\Rt$
are in any way related. However, we have:

\begin{prop} \label{special5prop}
  For all $j\in B\new$ there exists a unit $w_j\in\Rt^\times$ such
  that
  \[
       T_j \;\equiv\; w_j\cdot t_j \pmod{p}.
  \]
\end{prop}

\proof Let $R:=k[\epsilon]$ denote the ring of dual numbers. Fix some
$j_0\in B\new$ and let $\kappa:\Rt\to R$ be the unique $W(k)$-algebra
morphism which sends $t_{j_0}$ to $\epsilon$ and $t_j$ to $0$ for
$j\not=j_0$. Set $\tau_{j,R}:=\kappa(\tau_{j,\Rt})$. Then
$\tau_{j,R}=\tau_j+\epsilon\cdot \delta_j$ for an element $\delta_j\in k$. To
prove the proposition it suffices to show that
\begin{equation} \label{special5eq1}
   \delta_j \not=0  \quad\text{\rm if and only if}\quad j=j_0.
\end{equation}

Let $Y_R$ denote the equivariant deformation of $Y$ obtained from
pulling back the universal deformation $\Y$
along $\kappa$. The isomorphism class of $Y_R$ corresponds to a class
in $\HExt_G^1(\Ll_{Y/k},\Oo_Y)$. Via the isomorphism
\[
     \HExt_G^1(\Ll_{Y/k},\Oo_Y) \;\liso\; 
        \bigoplus_{j\in B\new}\; \T_{X,\tau_j}\otimes k(\tau_j),
\]
this class may be represented by a tuple $(\bar{\theta}_j)$, where
$\bar{\theta}_j$ is a tangent vector in $\tau_j$ (see the proof of
Proposition \ref{special3lem}). By the choice of the indeterminates
$t_j$ and the homomorphism $\kappa$ we have
\begin{equation} \label{special5eq2}
  \bar{\theta}_j  \not=0  \quad\text{\rm if and only if}\quad j=j_0.
\end{equation}
On the other hand, it follows from Proposition \ref{defdat3prop} that
the tuple $(\bar{\theta}_j)$, considered as a class in
\[
    H^1(X,\T_X(-\sum_{j\in B}\tau_j)) \;\cong\;
       H^0(X,\,\frac{\T_X(-\sum_{j\in B_0}\tau_j)}
                    {\T_X(-\sum_{j\in B}\tau_j)}\,) \;\cong\;
     \bigoplus_{j\in B\new}\;\T_{X,\tau_j}\otimes k(\tau_j),
\]
represents the isomorphism class of the deformation
$(X_R;\tau_{j,R})$. Therefore, 
\begin{equation} \label{special5eq3}
    \bar{\theta}_j \;=\; \delta_j\cdot\vf{x}|_{x=\tau_j}.
\end{equation}
Now \eqref{special5eq2} and \eqref{special5eq3} together imply
\eqref{special5eq1}. The proposition is proved.
\Endproof

To finish, let us explain briefly the motivation behind Proposition
\ref{special5prop}. Let $R$ be a complete discrete valuation ring of
mixed characteristic $(0,p)$. Let $k$ be the residue field of $R$
(which we assume algebraically closed) and $K$ its fraction field. Let
$(Z,V)$ be a special deformation datum over $X=\PP^1_k$ and $Y\to X$
the associated $G$-cover.  Furthermore, let $\tau_{j,R}\in R$ be
points on $X_R=\PP^1_R$ which lift the branch points $(\tau_j)$ of
$Y\to X$.  By the lifting property, there exists a unique lift $Y_R\to
X_R$ of $Y\to X$ with branch points $(\tau_{j,R})$.  Assuming that
$\zeta_p\in R$, the generic fiber $Y_K\to X_K=\PP^1_K$ is a tame
Galois cover with Galois group
\[
     G(K) \;\cong\; (\ZZ/p)^s\rtimes H.
\]
By construction, the cover $Y_K\to X_K$ has bad reduction: the special
fiber $Y$ is singular and the induced map $Y\to X$ is not separable.
However, after some blowing up we can find a certain nice model
$\tilde{Y}_R\to\tilde{X}_R$ over $R$ of $Y_K\to X_K$, called the
{\em stable model}, see \cite{Raynaud98} and \cite{bad}.

What can we say about the stable model, and how does it depend on the
choice of the branch points $\tau_{j,R}$? Let us say that the stable
reduction of $Y_K\to X_K$ is {\em nice} if the vanishing cycles of
$Y_R$ are resolved by blowing up $X_R=\PP^1_R$ in certain {\em
  disjoint} closed disks with center $\tau_{j,R}$, for $j\in
B-B\wild$ (in \cite{special} this property is called {\em special}).
If the stable reduction is nice, then the special fiber $\tilde{X}$ of
$\tilde{X}_R$ is a {\em comb}. More precisely, $\tilde{X}$ is a
semistable curve consisting of the central component $X$ and, for each
index $j\in B-B\new$, a tail $X_j$ meeting $X$ in $\tau_j$. Using
Proposition \ref{special5prop} one can show the following.

\begin{result}
  The stable reduction of $Y_K\to X_K$ is nice if and only if the
  branch points $\tau_{j,R}\in R$ are `sufficiently close' to the
  Teichm\"uller lift $[\tau_j]\in W(k)\subset R$.
\end{result}

In the case $s=\dim_{\FF_p}V=1$ the `if'-direction of this result was
proved in \cite{special}, using a very different kind of argument. In
\cite{bad} and still under the condition $s=1$, both directions of the
above result are proved, using Proposition \ref{special5prop}. The
case $s>1$ is similar but a bit more involved and will be dealt with
in a subsequent paper.

Roughly speaking, the singularities of $Y_R$ can be described in terms
of the image of the parameters $t_j$ in $R$ (under the classifying map
$\Rt\to R$ of the deformation $Y_R$). Therefore, Proposition
\ref{special5prop} provides a link between the singularities of $Y_R$
and the position of the branch points of $Y_R\to X_R$. It is somewhat
surprising that such a relation exists at all, because the dependence
of the cover $Y_R\to X_R$ on the branch points $\tau_{j,R}$ seems to
be of a more global nature.


\begin{appendix}

\section{Picard categories and Picard stacks} \label{pic}

In this first appendix we recall some basic facts about Picard
categories and Picard stacks. The main result we need is Proposition
\ref{pic2prop2}. References are \cite{SGA4}, Expos\'e XVIII and
\cite{Ulbrich84}.

\subsection{} \label{pic1}

A {\em (strictly commutative) Picard category} is a nonempty
monoid $\Pic$, together with a functor 
\[
     +:\;\Pic\times\Pic \;\To\; \Pic, \qquad (x,y) \;\longmapsto\; x+y
\]
and two functorial isomorphisms
\[
      \sigma:\; (x+y)+z \;\cong\; x+(y+z), \qquad
      \tau:\; x+y \;\cong\; y+x
\]
such that the following holds.
\begin{enumerate}
\item
    The isomorphisms $\sigma$ and $\tau$ make $+$ an {\em associative}
    and {\em strictly commutative} functor, in the sense of
    \cite{SGA4}, Expos\'e XVIII, \S 1.4.1.
\item
    For any object $y$ of $\Pic$, the functor $x\mapsto x+y$ is an
    equivalence of categories.
\end{enumerate}
Given two Picard categories $\Pic_1,\Pic_2$, an {\em additive functor}
from $\Pic_1$ to $\Pic_2$ is a functor $F: \Pic_1\to\Pic_2$ together
with a functorial isomorphism
\[
      F(x+y) \liso F(x)+F(y)
\]
which is compatible with the associativity and the commutativity laws,
see \cite{SGA4}, Expos\'e XVIII, \S 1.4.6. Given two additive functors
$F,G:\Pic_1\to\Pic_2$, a {\em morphism of additive functors} $u:F\to
G$ is a morphism of functors (automatically an isomorphism) such that
the diagram
\[\begin{CD}
        F(x+y)    @>{u_{x+y}}>> G(x+y)    \\
        @VVV                    @VVV      \\
        F(x)+F(y) @>{u_x+u_y}>> G(x)+G(y) \\
\end{CD}\]
commutes.  We denote by $\HOM(\Pic_1,\Pic_2)$ the corresponding
category of additive functors and by $\Hom(\Pic_1,\Pic_2)$ its set of
isomorphism classes. One can show that $\HOM(\Pic_1,\Pic_2)$ carries a
natural structure of Picard category.

\vspace{1ex} Let $A$ be a complex of abelian groups. We define a
Picard category $\PIC(A)$ as follows. Objects of $\PIC(A)$ are
$1$-cocycles, i.e.\ elements of $Z^1(A)=\Ker(A^1\pfeil{d}A^2)$. Given
two objects $x,y\in Z^1(A)$, the set of morphisms $\Hom(x,y)$ is the
set of elements $f\in A^0$ such that $d(f)=y-x$, modulo
$0$-coboundaries, i.e.\ elements of $B^0(A)=\Im(A^{-1}\pfeil{d}A^0)$.
The composition of two morphisms $f:x\to y$ and $g:y\to z$ is the sum
$f+g$. The functor $+$ is induced from the addition law of $A^1$. It
follows immediately from this definition that the group of
automorphisms of the `neutral object' of $\PIC(A)$ is identified with
$H^0(A)$, whereas the group of isomorphism classes of $\PIC(A)$ is
identified with $H^1(A)$. Note also that $\PIC(A)=\PIC(A^{[0,1]})$,
where $A^{[0,1]}$ denotes the complex of amplitude $[0,1]$ deduced
from $A$ such that $H^n(A^{[0,1]})=H^n(A)$ for $n=0,1$.

Given two complexes of abelian groups $A,B$, a homomorphism of complexes
$\varphi:A\to B$ gives rise to an additive functor
$\PIC(\varphi):\PIC(A)\to\PIC(B)$. The functor $\PIC(\varphi)$ is an
equivalence of categories if and only if $H^n(\varphi)$ is an
isomorphism for $n=0,1$. Given two homomorphisms $\varphi,\psi:A\to
B$, a homotopy $\varphi\sim\psi$ gives rise to an isomorphism of
additive functors $\PIC(\varphi)\cong\PIC(\psi)$. Therefore, the
association $A\mapsto \PIC(A)$ gives rise to a functor from the
derived category $\Der^{\geq 0}(\Ab)$ to the category of all (small)
Picard categories (morphisms in the latter category are isomorphism
classes of additive functors). It is shown in \cite{SGA4}, Expos\'e
XVIII, that this functor becomes an equivalence of categories when
restricted to the full subcategory $\Der^{[0,1]}(\Ab)$ of complexes of
amplitude $[0,1]$.

\begin{rem}
  Our definition of $\PIC(A)$ is a bit different from the definition
  used in \cite{SGA4}, Expos\'e XVIII. In {\em loc.cit.}, $\PIC(A)$ is
  only defined for a complex of amplitude $[-1,0]$. For the
  application of Picard categories in this paper, it seemed more
  convenient to to shift degrees by $1$ and to allow arbitrary
  complexes.
\end{rem}

\subsection{} \label{pic2}

Let $X$ be a topological space. (Actually, without changing
anything essential, we could let $X$ be an arbitrary site.)
We denote by $\Ab(X)$ the category of sheaves of abelian groups on
$X$. The total right derived functor of the global section functor
$\Gamma(X,\;\cdot\;)$ is denoted by $\RR\Gamma(X,\;\cdot\;)$. For
generalities about stacks, see \cite{Giraud} or \cite{LaumonMB}.

A {\em Picard stack} over $X$ is a stack $\PPic$ over $X$, together
with a morphism of $X$-stacks
\[
     +:\;\PPic\times_X\PPic \;\To\; \PPic, \qquad (x,y) 
            \;\longmapsto\; x+y
\]
and two functorial isomorphisms
\[
      \sigma:\; (x+y)+z \;\cong\; x+(y+z), \qquad
      \tau:\; x+y \;\cong\; y+x
\]
such that the following holds.
\begin{enumerate}
\item
    For each open subset $U\subset X$, the fiber $\PPic(U)$, together
    with the restrictions of $+$, $\sigma$ and $\tau$ to $U$, is a
    (strictly commutative) Picard category.
\item
    For each inclusion $U\subset V$ of open subsets, the restriction
    functor $\PPic(U)\to\PPic(V)$ is a morphism of Picard categories.
\end{enumerate}
Given two Picard stacks $\PPic_1,\PPic_2$, an additive functor from
$\PPic_1$ to $\PPic_2$ is an $X$-functor $F:\PPic_1\to\PPic_2$
together with functorial isomorphisms $F(x+y)\iso F(x)+F(y)$ whose
restriction to each fiber is an additive functor. Given two additive
functors $F,G:\PPic_1\to\PPic_2$, a morphism of additive functors is a
morphism of $X$-functors $u:F\to G$ (automatically an isomorphism)
whose restriction to all fibers is a morphism of additive functors.
We denote by $\HOM(\PPic_1,\PPic_2)$ the corresponding category of
additive functors and by $\Hom(\PPic_1,\PPic_2)$ its set of
isomorphism classes. It is easy to equip $\HOM(\PPic_1,\PPic_2)$ with
a natural structure of a Picard category. Moreover, one can show
that the $X$-groupoid
\[
      U \:\longmapsto\; \HOM(\PPic_1|_U,\PPic_2|_U)
\]
is itself a Picard stack, see \cite{SGA4}, Expos\'e XVIII. 

Let $\A$ be a complex of abelian sheaves on $X$. The association
\[
       U \;\longmapsto\; \PIC(\Gamma(U,\A))
\]
(where $U\subset X$ runs over all open subsets of $X$) gives rise to a
prestack $\PPic'(\A)$ over $X$. Let $\PPIC(\A)$ be the stack over $X$
associated to this prestack, see e.g.\ \cite{LaumonMB}, Lemme
(3.2). One checks that $\PPIC(\A)$ is a Picard stack, in a natural
way. For each open subset $U\subset X$, the natural functor
\begin{equation} \label{pic2eq1}
        \PIC(\Gamma(U,\A)) \;\To\; \PPIC(\A)
\end{equation}
is an additive functor. In general, it is not an isomorphism. 
  
A homomorphism $\varphi:\A\to\B$ of abelian sheaves gives rise to a
morphism $\PIC(\varphi):\PPIC(\A)\to\PPIC(\B)$ of Picard
stacks. Moreover, a homotopy $\varphi\sim\psi$ gives rise to a
isomorphism of additive functors
$\PIC(\varphi)\cong\PIC(\psi)$. Therefore, the association
$\A\mapsto\PPIC(\A)$ gives rise to a functor from the derived category
$\Der(X)$ to the category of all (small) Picard stacks on $X$
(morphisms in the latter category are isomorphism classes of additive
functors).

\begin{prop} \label{pic2prop2}
  Let $X$ be a topological space and $\A$ a sheaf of abelian groups on
  $X$ such that $H^n(\A)=0$ for $n<0$. Then we have a natural
  equivalence of Picard categories
  \[
       \PPIC(\A)(X) \;\cong\; \PIC(\RR\Gamma(X,\A)).
  \]
\end{prop}

This proposition seems to be well known. Since it is an important step
in the proof of Theorem \ref{extthm} and we could not find a suitable
reference, we give a proof. 

\proof 
The hyper-cohomology spectral sequence and the assumption $H^n(\A)=0$
for $n<0$ show that 
\[
       \HH^n(X,\A^{[0,1]}) \;=\; \HH^n(X,\A) \qquad
           \text{for $n=0,1$.} 
\]
We may therefore assume that $\A=\A^{[0,1]}$. Let $\U=(U_i)_{i\in I}$
be an open covering of $X$. We choose a well-ordering on the index set
$I$. Let $K_{\U}:=C^\bullet(\U,\A)$ be the double complex whose $n$th
column (for $n=0,1$) is the $\check{\rm C}$ech cochain complex of
$\A^n$ with respect to $\U$:
\[
       K_{\U}:\quad\left\{\quad
  \begin{array}{ccc} 
    \Prod_i\;\Gamma(U_i,\A^0) & \lpfeil{d} & 
             \Prod_i\;\Gamma(U_i,\A^1) \\ 
    \rupfeil{\partial}  & &  \rupfeil{\partial} \\ 
     \Prod_{i<j}\;\Gamma(U_{i,j},\A^0) & \lpfeil{-d} & 
             \Prod_{i<j}\;\Gamma(U_{i,j},\A^1) \\
    \rupfeil{\partial}  & &  \rupfeil{\partial} \\ 
     \Prod_{i<j<k}\;\Gamma(U_{i,j,k},\A^0) & \lpfeil{d} &  
              \Prod_{i<j<k}\;\Gamma(U_{i,j,k},\A^1) \\
    \rupfeil{\partial}  & &  \rupfeil{\partial} \\
      \cdots  & & \cdots 
  \end{array}
    \right.
\]
We define a morphism of Picard categories
\begin{equation} \label{pic2eq4}
   \PIC(\Tot(K_{\U})) \;\liso\; \PPIC(\A)(X),
\end{equation}
as follows. An object of $\PIC(\Tot(K_{\U}))$ is a datum $(f_i;g_{i,j})$,
with $f_i\in\Gamma(U_i,Z^1(\A))$ and $g_{i,j}\in\Gamma(U_{i,j},\A^0)$,
such that 
\begin{equation} \label{pic2eq5}
   d(g_{i,j}) \;=\; f_j|_{U_{i,j}}-f_i|_{U_{i,j}}
\end{equation}
for all $i<j$ and
\begin{equation} \label{pic2eq6}
   g_{i,j} - g_{i,k} + g_{j,k} \;=\; 0
\end{equation}
for all $i<j<k$. Let $\tilde{f}_i$ denote the object of
$\PPIC(\A)(U_i)$ corresponding to $f_i$. By \eqref{pic2eq5}, $g_{i,j}$
corresponds to an isomorphism
$\tilde{g}_{i,j}:\tilde{f}_i|_{U_{i,j}}\iso\tilde{f}_j|_{U_{i,j}}$.
Now \eqref{pic2eq6} means that these isomorphisms satisfy the cocycle
relation $\tilde{g}_{j,k}\circ\tilde{g}_{i,j}=\tilde{g}_{i,k}$. In
other words, $(\tilde{f}_i;\tilde{g}_{i,j})$ is a patching datum with
values in
$\PPIC(\A)$. Since $\PPIC(\A)$ is a stack, there exists an object
$\tilde{f}$ of $\PPIC(\A)(X)$ together with isomorphisms
$\alpha_i:\tilde{f}|_{U_i}\iso\tilde{f}_i$ such that
$\tilde{g}_{i,j}=\alpha_j\circ\alpha_i^{-1}$. By definition,
$\tilde{f}$ is the image of $(f_i,g_{i,j})$ under 
\eqref{pic2eq4}. 

Let $(f_i',g_{i,j}')$ be another object of $\PIC(\Tot(K_{\U}))$, and let
$\tilde{f}'$ be the corresponding object of $\PPIC(\A)(X)$. A
homomorphism from $(f_i,g_{i,j})$ to $(f_i',g_{i,j}')$ is a datum
$(h_i)$, with $h_i\in\Gamma(U_i,\A^0)$, such that
\begin{equation} \label{pic2eq7}
   d(h_i) \;=\; f_i' - f_i
\end{equation}
for all $i$ and 
\begin{equation} \label{pic2eq8}
  h_j|_{U_{i,j}} - h_i|_{U_{i,j}} \;=\; g_{i,j}' - g_{i,j}
\end{equation}
for all $i<j$. Equation \eqref{pic2eq7} shows that $h_i$ corresponds to
an isomorphism $\tilde{h}_i:\tilde{f}_i\iso\tilde{f}_j$. Moreover, by
\eqref{pic2eq8} the diagram
\[\begin{CD}
   \tilde{f}_i|_{U_{i,j}} @>{\tilde{h}_i}>> \tilde{f}'_i|_{U_{i,j}} \\
   @V{\tilde{g}_{i,j}}VV                    @V{\tilde{g}'_{i,j}}VV  \\
   \tilde{f}_j|_{U_{i,j}} @>{\tilde{h}_j}>> \tilde{f}'_j|_{U_{i,j}} \\
\end{CD}\]
commutes for all $i<j$. In other words, $(\tilde{h}_i)$ is an
isomorphism of patching data. Using again that $\PPIC(\A)$ is a stack,
we see that there exists a unique isomorphism
$\tilde{h}:\tilde{f}\iso\tilde{f}'$ such that
$\tilde{h}_i\circ\alpha_i=\alpha_i'\circ\tilde{h}|_{U_i}$. By
definition, $\tilde{h}$ is the image of $(h_i)$ under
\eqref{pic2eq4}. This finishes the definition of \eqref{pic2eq4} as a
functor. We leave it to the reader to check that \eqref{pic2eq4} is
indeed a morphism of Picard categories.

Clearly, the definition of \eqref{pic2eq4} is compatible with taking
refinements of the covering $\U$. Therefore, we obtain a morphism of
Picard categories
\begin{equation} \label{pic2eq9}
  \dirlim{\U} \PIC(\Tot(K_{\U})) \;\To\;
      \PPIC(\A)(X).
\end{equation}
We claim that \eqref{pic2eq9} is an isomorphism. Indeed, the
discussion of the previous paragraph, leading to the definition of
\eqref{pic2eq4}, shows that $\PIC(\Tot(K_{\U}))$ is isomorphic to the
category of patching data for the covering $\U$, with values in the
prestack $\PPIC(\A)'$. On the other hand, the category $\PPIC(\A)(X)$
is the direct limit over the categories of such patching data, where
the limit is taken over all possible coverings $\U$; this follows from
the construction of a stack associated to a prestack, see e.g.\ 
\cite{LaumonMB}, \S 3. This proves the claim.

To finish the proof of the proposition, it suffices to show that the
natural morphisms
\[
     \Tot(K_{\U}) \;\To\; \RR\Gamma(X,\A) 
\]
induces isomorphisms on cohomology
\[
      \dirlim{\U} H^n(\Tot(K_{\U})) \;\liso\; \HH^n(X,\A) 
\] 
for $n=0,1$. This is proved in two steps. First, one compares the two
spectral sequences which compute the cohomology of $\Tot(K_{\U})$ on
the one hand and the hyper-cohomology groups $\HH^n(X,\A)$ on the other
hand. Then one uses the well known fact that $\check{\rm
  C}$ech-cohomology agrees with ordinary sheaf cohomology in degree
$n=0,1$ (see e.g.\ \cite{Hartshorne}, Ex.\ III.4.4). We omit the
details.  \Endproof


\section{Group cohomology for affine flat group schemes} \label{GRmod}

We show how to compute the cohomology of an affine group flat group
scheme in terms of cocycles and coboundaries, just as for abstract
groups. Reference is \cite{SGA3}, Expos\'e I. 

\subsection{} \label{GRmod1}

Throughout this section, we fix a commutative ring and an affine flat
$R$-group scheme $G=\Spec \Oo_G$. We denote by
$\Delta:\Oo_G\to\Oo_G\otimes_R\Oo_G$ the comultiplication and by
$e:\Oo_G\to R$ the counit of $G$.

A {\em (right) $G$-$R$-module} is an $R$-module $M$ together with an
$R$-linear morphism $\mu_M:M\to\Oo_G\otimes_R M$ such that
$(\Id_{\Oo_G}\otimes\mu_M)\circ\mu_M=(\Delta\otimes\Id_M)\circ\mu_M$
and $(e\otimes\Id_M)\circ\mu_M=\Id_M$. For each $R$-algebra $R'$ and
$\sigma\in G(R')$ we obtain an $R'$-linear automorphism $m\mapsto
m^\sigma$ of $M':=M\otimes_R R'$ such that
$m^{\sigma\tau}=(m^\sigma)^\tau$.  We shall denote by $\Mod(R,G)$ the
category of $G$-$R$-modules, by $\Ko^+(R,G)$ the (triangulated)
category of bounded below cochain complex in $\Mod(R,G)$ (modulo
homotopy) and by $\Der^+(R,G)$ the derived category of $\Ko^+(R,G)$.

Given a $G$-$R$-module $M$, the {\em invariant} $R$-submodule $M^G$ is
the set of all $m\in M$ such that $\mu_M(m)=1\otimes m$, or, what is
equivalent, $m^\sigma=m$ for all $R'$ and $\sigma\in G(R')$.  The
functor $M\mapsto M^G$ from $\Mod(R,G)$ to the category of abelian
groups is obviously additive and left exact.  We denote its $n$th
right derived functor by $H^n(G,\,\cdot\,)$ and its total right
derived functor by $\RR^G$. (For the existence of enough injectives in
$\Mod(G,R)$, see the proof of Lemma \ref{coindlem} below.)

\subsection{} \label{GRmod2}

Given an $R$-module $M$, we set
\[
     \tilde{M} \;:=\; \Oo_G\otimes_R M.
\]
The map $\Delta\otimes\Id_M:\tilde{M}\to\Oo_G\otimes_R\tilde{M}$
gives $\tilde{M}$ the structure of a (right) $G$-$R$-module. A
$G$-$R$-module which is isomorphic to $\tilde{M}$ for some $R$-module
$M$ is called {\em coinduced}. If $G$ is a finite group then this
agrees with the usual definition of coinduced modules.

Let $M$ be an $R$-module $M$, $P$ a $G$-$R$-module and $\varphi:P\to
M$ an $R$-linear morphism. Then
$\tilde{\varphi}:=(\Id_{\Oo_G}\otimes\varphi)\circ\mu_P:P\to\tilde{M}$
is easily checked to be $G$-equivariant. One checks that this
construction yields a natural isomorphism
\begin{equation} \label{coindeq1}
   \Hom_R(P,M) \liso \Hom_G(P,\tilde{M}).
\end{equation}
(The inverse of \eqref{coindeq1} is defined as follows: given a
$G$-equivariant homomorphism $\psi:P\to\tilde{M}$,
$\varphi:=(e\otimes\Id_M)\circ\psi:P\to M$ is an $R$-linear morphism
such that $\psi=\tilde{\varphi}$.) Moreover, the isomorphism
\eqref{coindeq1} makes the functor $M\mapsto\tilde{M}$ a right adjoint
of the forgetful map from $\Mod(R,G)$ to $\Mod(R)$.

\begin{lem} \label{coindlem}
  For any $R$-module $M$ we have
  \[
      H^n(G,\tilde{M}) \;=\; 
       \begin{cases}
          \;M & \;\text{for $n=0$},\\
          \;\;0 & \;\text{for $n>0$}.\\
       \end{cases}
  \]
\end{lem}

\proof For $n=0$, the lemma is equivalent to the exactness of the
sequence
\begin{equation} \label{coindeq2}
  0 \;\to\; M \;\To\; \tilde{M} \;\To\; \Oo_G\otimes_R\tilde{M}
\end{equation}
(the first arrow sends $m$ to $1\otimes m$ and the second $a\otimes m$
to $\Delta(a)\otimes m-1\otimes a\otimes m$). Now \eqref{coindeq2} is
exact on the left because $R\to\Oo_G$ is flat, by assumption.
Exactness in the middle is proved using the properties of the counit
$e:\Oo_G\to R$. Hence the lemma holds for $n=0$.

Choose an injective resolution $M\to I^0\to I^1\to\cdots$ of the
$R$-module $M$. The
functor $M\mapsto \tilde{M}$, being the right adjoint of an exact
functor, preserves injectives, see e.g.\ \cite{Weibel}, Proposition
2.3.11. Therefore, $\tilde{M}\to\tilde{I}^0\to\tilde{I}^1\to\cdots$ is
an injective resolution of the $G$-$R$-module $\tilde{M}$. Now the
general case of the lemma follows from the case $n=0$.
\Endproof

\subsection{} \label{GRmod3}

Let $M$ be a $G$-$R$-module. In order to compute the cohomology groups
$H^n(G,M)$, it suffices to write down a resolution of $M$ by coinduced
$G$-$R$-modules. We do this as follows. Set $B^{-1}(G,M):=M$ and
define inductively $B^n(G,M):=(B^{n-1}(M))^\sim$ for all $n\geq 0$.
As $R$-modules, we simply get 
\[
  B^n(G,M) \;=\;  
     \underbrace{\Oo_G\otimes_R\cdots\otimes_R\Oo_G}_
        {(n+1)\times}\otimes_R M.
\]
We define differentials $\partial:B^n(G,M)\to
B^{n+1}(G,M)$  by setting 
\begin{align*}
  \partial(a_0\otimes\cdots\otimes a_n\otimes m) \;:=\; &
    \sum_{\nu=0}^n\;(-1)^\nu\;a_0\otimes\cdots\otimes\Delta(a_\nu)
       \otimes\cdots\otimes a_n\otimes m  \\
        & \;+\; (-1)^{n+1}\;a_0\otimes\cdots\otimes a_n\otimes\mu_M(m).
\end{align*}
Note that $\partial:M=B^{-1}(G,M)\to\tilde{M}=B^0(G,M)$ equals
$\mu_M$. 
It is easy to check that 
\begin{equation} \label{canres}
    M \;\lpfeil{\mu_M}\; B^0(G,M) \;\lpfeil{\partial}\; B^1(G,M) 
      \;\lpfeil{\partial}\;\cdots
\end{equation}
is an (augmented) complex a $G$-$R$-modules. (In fact, \eqref{canres}
is the (augmented) cochain complex associated to $M$ and the pair of
adjoint functors $\Mod(G,R)\rightleftarrows\Mod(R)$, see
\cite{Weibel}, \S 8.6.) Moreover, \eqref{canres} is exact. This
follows immediately from the existence of the homotopy 
\[
  s:\;\left\{\;
    \begin{array}{ccc}
      B^{n+1}(G,M) & \;\To\; & B^n(G,M) \\
      a_0\otimes\cdots\otimes a_n\otimes m & \;\longmapsto\; &
            e(a_0)a_1\otimes\cdots\otimes a_n\otimes m
    \end{array}\right..
\]

Applying the functor $M\mapsto M^G$ to the resolution \eqref{canres}
defines an augmented complex of $R$-modules
\begin{equation} \label{cochains}
    M^G \;\To\; C^0(G,M) \;\lpfeil{\partial}\; C^1(G,M) 
      \;\lpfeil{\partial}\;\cdots
\end{equation}
Elements of $C^n(G,M)$ are called {\em $n$-cochains} with values in
$M$. Likewise, elements of $Z^n(G,M):=\Ker(\partial)$ (resp.\ of
$B^n(G,M):=\Im(\partial)$) are called {\em cocycles} (resp.\ {\em
  coboundaries}). Note that we have an isomorphism of $R$-modules
\[
  C^n(G,M) \;\cong\;  
    \underbrace{\Oo_G\otimes_R\cdots\otimes_R\Oo_G}_{n\times}\otimes_R M,
\]
such that the canonical injection $C^n(G,M)\inj B^n(G,M)$ sends the
element $a_1\otimes\cdots\otimes a_n\otimes m$ to the element
$1\otimes a_1\otimes\cdots\otimes a_n\otimes m$. Given an $R$-algebra
$R'$, an $n$-cochain $\varphi\in C^n(G,M)$ gives rise to a function
\[
   G(R')\times\cdots\times G(R') \To M'=M\otimes_R R', \qquad
   \underline{\sigma}=(\sigma_1,\ldots,\sigma_n) \longmapsto 
     \varphi_{\underline{\sigma}}.
\]
Now $\varphi$ is a cocycles (resp.\ a coboundary) if and only if this
function is a cocycle (resp.\ a coboundary) in the traditional sense,
for all $R$-algebras $R'$ (again, it suffices to take $R'$ flat over
$R$). For instance, a $1$-cochain $\varphi$ is a cocycle if and only
if 
\[
    \varphi_{\sigma,\tau} \;=\; \varphi_{\sigma}^\tau + \varphi_\tau
\]
holds for all $\sigma,\tau\in G(R')$. It is a coboundary if and only
if for all $R'$ there exists an element $m\in M'$ such that
$\varphi_\sigma=m^\sigma-m$ holds for all $\sigma\in G(R')$. 

It follows from Lemma \ref{coindlem} that 
\begin{equation} \label{cochains2}
        H^n(G,M) \;=\; H^n(C^\bullet(G,M))
\end{equation}
for all $G$-$R$-modules $M$ and all $n\geq 0$. The next proposition is a
slight generalization of \eqref{cochains2}.

\begin{prop} \label{GRmod3prop}
  Let $M^\bullet\in \Ko^+(R,G)$ be a bounded below complex of
  $G$-$R$-modules. Then we have a natural isomorphism of derived
  complexes 
  \[
      \RR^G(M^\bullet) \;\cong\; \Tot(C^\bullet(G,M^\bullet)).
  \]
\end{prop}

\proof Let $K$ denote the double complex $B^\bullet(G,M^\bullet)$. The
$q$th row of $K$ is exact except at degree $p=0$, where the cohomology
is $M^q$. Therefore, the spectral sequence associated to
$K$ (filtered by rows) shows that the augmentation $M^\bullet\to K$
gives rise to a quasi-isomorphism
\[
         M^\bullet \To \Tot(K).
\]
By definition, the complex $\Tot(K)$ consists entirely of coinduced
$G$-$R$-modules, which are acyclic with respect to taking
$G$-invariants, by Lemma \ref{coindlem}. Therefore, \cite{Weibel},
Theorem 10.5.9 implies
\[
    \RR^G(M^\bullet) \;=\; \Tot(K)^G \;=\;
        \Tot(C^\bullet(G,M^\bullet)).
\]
This finishes the proof of the proposition.
\Endproof


\section{Sheaves of $G$-$\Oo_Y$-modules} \label{GOY}

The goal of this last appendix is to review the definition of
equivariant hyperext groups and the construction of the two spectral
sequences \eqref{GOY4eq4} and \eqref{GOY5eq3}. The standard reference
is \cite{GroTohoku}.

\subsection{} \label{GOY2}

Let $G\to S$ and $Y\to S$ be as in \S \ref{constr}. Let
$\lambda:G\times_S Y\to Y$ (resp.\ $p:G\times_S Y\to Y$) denote the
morphism defining the action of $G$ on $Y$ (resp.\ the second
projection). Given a sheaf of $\Oo_Y$-modules $\F$, a {\em lift} of the
$G$-action from $Y$ to $\F$ is given by an isomorphism
$\lambda^*\F\to p^*\F$ which satisfies certain obvious axioms, see
e.g.\ \cite{MumfordAV}, \S III.12.

Let $\F$ and $\G$ be $G$-$\Oo_Y$-modules. Let $\Hom_Y(\F,\G)$ denote
the $R$-module of $\Oo_Y$-linear (but not necessarily $G$-equivariant)
homomorphisms from $\F$ to $\G$. It carries a natural structure of
$G$-$R$-module, defined as follows. Let $R'$ be a flat $R$-algebra and
$\sigma\in G(R')$. Since $R'$ is flat over $R$ we have a natural
isomorphism
\[
       \Hom_Y(\F,\G)\otimes_R R' \;\cong\; \Hom_{Y'}(\F',\G').
\]
Given $f:\F'\to\G'\in\Hom_{Y'}(\F',\G')$, we define $f^\sigma$ via
the following commutative diagram:
\[\begin{CD}
      \F'         @>{f^\sigma}>>   \G'           \\
      @V{\varphi_\sigma}VV  @VV{\varphi_\sigma}V \\
      \sigma^*\F' @>>{\sigma^*f}>  \sigma^*\G'   \\
\end{CD}\]
Note that an $\Oo_Y$-linear morphism $f:\F\to\G$ is $G$-equivariant if
and only if it is invariant under the $G$-action just defined, i.e.\ 
\[
       \Hom_G(\F,\G) \;=\; \Hom_Y(\F,\G)^G.
\]

\subsection{} \label{GOY3}

Let $\F$ be an $\Oo_Y$-module. We set $\tilde{\F}:=\mu_*p^*\F$ and
claim that $\tilde{\F}$ carries a natural structure of
$G$-$\Oo_Y$-module.  To define a $G$-action on $\tilde{\F}$, let
$S'=\Spec R'$ be an affine $S$-scheme and $\sigma\in G(S')$. As usual,
we indicate base change to $S'$ with a prime (e.g.\ $Y'=Y\times_S S'$)
and identify $\sigma$ with the automorphism of $Y'$ induced from
$\sigma$ via the action of $G$ on $Y$. Also, we let
$t_\sigma:G'\times_{S'} Y'\iso G'\times_{S'} Y'$ denote the
automorphism induced from left translation by $\sigma$ on the first
factor. Since $p'\circ t_\sigma=p'$, we have a natural isomorphism
\begin{equation} \label{GOY3eq1}
  {p'}^*\F' \liso t_\sigma^*{p'}^*\F'.
\end{equation}
Using the commutative diagram
\begin{equation} \label{GOY3eq2}
\begin{CD}
  G'\times_{S'}Y'   @>{t_\sigma}>>  G'\times_{S'}Y' \\
  @V{\mu'}VV                        @VV{\mu'}V      \\
  Y'                @>{\sigma}>>    Y'              \\
\end{CD}
\end{equation}
we obtain an isomorphism
$\varphi_\sigma:\tilde{\F}'\iso\sigma^*\tilde{\F}'$ as follows:
\begin{equation} \label{GOY3eq3}
  \tilde{\F}' \;=\; \mu_*'{p'}^*\F' \liso \mu_*'t_\sigma^*{p'}^*\F'
       \liso \sigma^*\mu_*'{p'}^*\F' \;=\; \sigma^*\tilde{\F}'.
\end{equation}
One checks that the definition of $\varphi_\sigma$ is functorial in
$S'$ and satisfies the rule
$\varphi_{\sigma\tau}=(\tau^*\varphi_\sigma)\circ\varphi_\tau$. This
defines a lift of the $G$-action from $Y$ to $\tilde{\F}$. we may
therefore regard $\tilde{\F}$ as a $G$-$\Oo_Y$-module.  A
$G$-$\Oo_Y$-module isomorphic to $\tilde{\F}$ for some $\Oo_Y$-module
$\F$ is called {\em coinduced}. Compare with \S \ref{GRmod2}.

\begin{prop} \label{GOY3prop}
\begin{enumerate}
\item
  Given an $\Oo_Y$-module $\F$ and a $G$-$\Oo_Y$-module $\G$, we have a 
  natural isomorphism of $G$-$R$-modules
  \[
       \Hom_Y(\G,\tilde{\F}) \;\cong\; \widetilde{\Hom_Y(\G,\F)}
  \]
  (the right hand side is the $G$-$R$-module coinduced from the
  $R$-module $\Hom_Y(\G,\F)$).
\item 
  The functor $\F\longmapsto\tilde{\F}$ is the right adjoint of
  the forgetful functor $\Mod(Y,G)\to\Mod(Y)$.
\end{enumerate}
\end{prop}

\proof 
Given a $G$-$\Oo_Y$-module $\G$, we have natural isomorphisms of
$R$-modules
\begin{equation} \label{GOY3eq4}
  \Hom_Y(\G,\F) \;\cong\; \Hom_{G\times Y}(\mu^*\G,p^*\F) \;\cong\;
    \Hom_{G\times Y}(p^*\G,p^*\F) \;\cong\; \widetilde{\Hom_Y(\G,\F)}.
\end{equation}
Here the first isomorphism is obtained from the adjointness of $\mu_*$
and $\mu^*$, the second isomorphism from the $G$-action on $\G$ and
the third isomorphism exists because $G\to S$ is flat. Note that both
the first and the last term in \eqref{GOY3eq4} carry a natural
structure of $G$-$R$-module. One checks that \eqref{GOY3eq4} defines
an isomorphism of $G$-$R$-modules. This proves (i). Taking
$G$-invariants and using Lemma \ref{coindlem} for $n=0$, we obtain an
isomorphism of $R$-modules
\begin{equation} \label{GOY3eq5}
   \Hom_G(\G,\tilde{\F}) \;\cong\; \Hom_Y(\G,\F).
\end{equation}
A tedious but elementary verification shows that this isomorphism
makes $\F\mapsto\tilde{\F}$ the right adjoint of the forgetful functor
$\Mod(Y,G)\to\Mod(Y)$. 
\Endproof

\begin{cor} \label{GOY3cor}
\begin{enumerate}
\item
  The category $\Mod(Y,G)$ has enough injectives.
\item
  Given a $G$-$\Oo_Y$-module $\G$, the functor 
  \[
      \Hom_Y(\G,\;\cdot\;):\; \Mod(Y,G) \;\To\; \Mod(R,G)
  \]
  sends injective $G$-$\Oo_Y$-modules to $G$-$R$-modules which are
  acyclic with respect to the functor $M\mapsto M^G$.
\end{enumerate}
\end{cor}

\subsection{} \label{GOY4}

Let $\A,\B$ be complexes of $G$-$\Oo_Y$-modules. We assume that $\A$ is
bounded below and that $\B$ is bounded in both directions. Then the
{\em total Hom complex} $\Hom_G^\bullet(\A,\B)$ is
a (bounded above) complex of abelian groups whose $n$th cohomology
group is isomorphic to the group of homomorphisms $\A\to\B[-n]$ up to
homotopy, i.e.
\[
       H^n(\Hom_G^\bullet(\A,\B)) \;=\; \Hom_{\Ko^+(Y,G)}(\A,\B[-n]).
\]
See e.g.\ \cite{Weibel}, 2.7.4.  Let $\Ko^-(\Ab)$ denote the category
of bounded above complexes of abelian groups, up to homotopy, and
$\Der^-(\Ab)$ its derived category.  The functor
$\Hom_G^\bullet(\A,\,\cdot\,):\Ko^b(G,Y)\to\Ko^-(\Ab)$ is a morphism
of triangulated categories and has a total right derived functor,
which we denote by $\RHom_G(\A,\,\cdot\,):\Der^b(Y,G)\to\Der^-(\Ab)$.
The $n$th {\em hyperext} of $\A$ and $\B$ is defined as
\[
     \HExt_G^n(\A,\B) \;:=\; H^n(\RHom_G(\A,\B))
\]
see e.g.\ \cite{Weibel}, \S 10.7. It follows from standard arguments
that the functor $\RHom_G(\,\cdot\,,\B):\Der^+(Y,G)\to\Der^-(\Ab)$ is a
morphism of triangulated categories.

Ignoring the $G$-action, we may as well define the total Hom complex
$\Hom_Y^\bullet(\A,\B)$. By \S \ref{GOY2}, the terms of
$\Hom_Y^\bullet(\A,\B)$ carry a natural structure of $G$-$R$-modules
such that 
\begin{equation} \label{GOY4eq2}
   \Hom_G^\bullet(\A,\B) \;=\; \Hom_Y^\bullet(\A,\B)^G.
\end{equation}
This formula displays the functor $\Hom_G^\bullet(\A,\;\cdot\;)$ as
the composition of two morphisms of triangulated categories
\[
    \Ko^b(Y,G) \;\To\; \Ko^-(R,G) \;\To\; \Ko^-(\Ab).
\]
It follows from Corollary \ref{GOY3cor} and \cite{Weibel}, Theorem
10.8.2 that 
\begin{equation} \label{GOY4eq3}
    \RHom_G(\A,\B) \;\cong\; \RR^G(\RHom_Y(\A,\B)).
\end{equation}
Therefore, the hyperext group $\HExt_G^n(\A,\B)$ can be computed via
the Grothendieck spectral sequence
\begin{equation} \label{GOY4eq4}
    E_2^{p,q} := H^p(G,\HExt_Y^q(\A,\B)) \;\Longrightarrow\;
       \HExt_G^{p+q}(\A,\B).
\end{equation}

\subsection{} \label{GOY5}

Now assume that $G\to S$ is finite and that $Y$ can be covered by
affine $G$-stable open subsets. Then by \cite{MumfordAV}, Theorem
12.1, the quotient scheme $X:=Y/G$ exists and the natural projection
$\pi:Y\to X$ is finite. Given a $G$-$\Oo_Y$-module $\F$, one defines a
sheaf of $\Oo_X$-modules $\F^G$ such that
\[
       \Gamma(V,\F^G) \;=\; \Gamma(\pi^{-1}(V),\F)^G.
\]
Let $\F$, $\G$ be two $G$-$\Oo_Y$-modules. In view of \S
\ref{GOY2} it is clear that the sheaf $\HOm_Y(\F,\G)$ is endowed with a
natural action of $G$, i.e.\ with a structure of $G$-$\Oo_Y$-module. We
set
\[
     \HOm_G(\F,\G) \;:=\; \HOm_Y(\F,\G)^G.
\]
By definition, we have
\begin{equation} \label{GOY5eq1}
    \Hom_G(\F,\G) \;=\; \Gamma(X,\HOm_G(\F,\G)).
\end{equation}

\begin{lem} \label{GOY5lem}
  The additive functor $\HOm_G(\F,\,\cdot\,)$ sends injective
  $G$-$\Oo_Y$-modules to $\Gamma(X,\,\cdot\,)$-acyclic $\Oo_X$-modules.
\end{lem}

\proof
By the construction of injective objects of $\Mod(Y,G)$ in \S
\ref{GOY3}, it suffices to proof the lemma for 
$G$-$\Oo_Y$-modules of the form $\tilde{\I}$, where $\I$ is an
injective $\Oo_Y$-module. Using Proposition \ref{GOY3prop} (ii) one
shows that
\[
   \HOm_G(\F,\tilde{\I}) \;=\; \pi_*\HOm_Y(\F,\I).
\]
It is well known that $\HOm_Y(\F,\,\cdot\,)$ sends injective to
$\Gamma(Y,\,\cdot\,)$-acyclic $\Oo_Y$-modules. But since $\pi$ is
finite the functor $\pi_*$ is exact and so $\pi_*\HOm_Y(\F,\I)$ is
$\Gamma(X,\,\cdot\,)$-acyclic. This proves the lemma.  \Endproof

Given two complexes of $G$-$\Oo_Y$-modules $\A$ and $\B$ (with $\A$
bounded below and $\B$ bounded), one defines the (bounded above)
complex of $\Oo_X$-modules $\HOm_G^\bullet(\A,\B)$. We let
$\RHOm_G(\A,\B)$ denote the total right derived functor of
$\HOm_G^\bullet(\A,\,\cdot\,)$, evaluated at $\B$, and set
\[
     \EXt_G^n(\A,\B) \;:=\; H^n(\RHOm_G(\A,\B)).
\]
As in the previous subsection, it follows from \eqref{GOY5eq1} and
Lemma \ref{GOY5lem} that
\begin{equation} \label{GOY5eq2}
   \RHom_G(\A,\B) \;\cong\; \RR\Gamma(X,\RHOm_G(\A,\B)).
\end{equation}
In particular, there exists a spectral sequence
\begin{equation} \label{GOY5eq3}
   H^p(X,\EXt_G^q(\A,\B)) \;\Longrightarrow\; \HExt_G^{p+q}(\A,\B).
\end{equation}
Following \cite{GroTohoku}, we call \eqref{GOY5eq3} the {\em
  local-global spectral sequence} for $\HExt_G^*$.

\end{appendix}



\vspace{5mm}
\flushright{Mathematisches Institut\\
            Universit\"at Bonn\\
            Beringstr. 1\\
            53115 Bonn, Germany\\
            wewers@math.uni-bonn.de}

\end{document}